\documentclass[english,10pt]{article}

\usepackage{xcolor}
\usepackage{amsmath}
\usepackage{graphics}
\usepackage{epsfig}
\usepackage{geometry}
\usepackage{pdfpages}
\usepackage{amssymb}
\usepackage[latin1]{inputenc}
\usepackage{babel,indentfirst,amsfonts,array}

\newtheorem{thm}{Theorem}[section]
\newtheorem{nthm}{Numerical Evidence}[section]
\newtheorem{lemme}[thm]{Lemma}

\newtheorem{defi}[thm]{Definition}

\def\P{\mathbb P}
\def\E{\mathbb E}
\def\T{\mathbb T}
\def\R{\mathbb R}
\def\C{\mathbb C}
\def\Q{\mathbb Q}
\def\Z{\mathbb Z}
\def\N{\mathbb N}
\def\({\left(}
\def\){\right)}
\def\[{\left[}
\def\]{\right]}
\def\fin{\hfill\square}

\def\fin{\hfill $\square$}

\let\oldsqrt\sqrt
\def\DHLhksqrt#1#2{%
\setbox0=\hbox{$#1\oldsqrt{#2\,}$}\dimen0=\ht0
\advance\dimen0-0.2\ht0
\setbox2=\hbox{\vrule height\ht0 depth -\dimen0}%
{\box0\lower0.4pt\box2}}
\title{\textsc{Self-similar measures and the Rajchman property}
\author{Julien Br\'emont}
\date{Universit\'e Paris-Est-Cr\'eteil,~novembre 2020}
}

\begin{document}
\maketitle
\setcounter{page}{1}

\begin{abstract}
For Bernoulli convolutions, the convergence to zero of the Fourier transform at infinity was characterized by successive works of Erd\"os \cite{erdos} and Salem \cite{salem}. We provide a quasi-complete extension of these results to general self-similar measures on the real line.\end{abstract}

\footnote{
\begin{tabular}{l}\textit{AMS $2010$ subject classifications~: 11K16, 37A45, 42A38, 42A61, 60K20.} \\
\textit{Key words and phrases~: Rajchman measure, self-similar measure, Pisot number, Plastic number.} 
\end{tabular}}

\section{Introduction}
{\it Rajchman measures.} In the present article we consider the question of extending some classical results on Bernoulli convolutions to a more general context of self-similar measures. For a Borel probability measure $\mu$ on $\R$, define its Fourier transform as~:

$$\hat{\mu}(t)=\int_{\R}e^{2i\pi tx}~d\mu(x),~t\in\R.$$

\smallskip
\noindent
We say that $\mu$ is Rajchman, whenever $\hat{\mu}(t)\rightarrow0$, as $t\rightarrow+\infty$. When $\mu$ is a Borel probability measure on the torus $\T=\R\backslash\Z$, we introduce its Fourier coefficients, defined as~:

$$\hat{\mu}(n)=\int_{\T}e^{2i\pi nx}~d\mu(x),~n\in\Z.$$

\smallskip
\noindent
In this study, starting from a Borel probability measure $\mu$ on $\R$, Borel probability measures on $\T$ will naturally appear, quantifying the non-Rajchman character of $\mu$. 

\medskip
\noindent
For a Borel probability measure $\mu$ on $\R$, the Rajchman property holds for example if $\mu$ is absolutely continuous with respect to Lebesgue measure ${\cal L}_{\R}$, by the Riemann-Lebesgue lemma. The situation can be more subtle and for instance there exist Cantor sets of zero Lebesgue measure and even of zero-Hausdorff dimension which support a Rajchman measure; cf Menshov \cite{menshov}, Bluhm \cite{bluhm}. Questions on the Rajchman property of a measure naturally arise in Harmonic Analysis, for example when studying sets of multiplicity for trigonometric series; cf Lyons \cite{seven} or Zygmund \cite{zyg}. We shall say a word on this topic at the end of the article. A classical counter-example is the uniform measure $\mu$ on the standard middle-third Cantor set, which is a continuous singular measure, not Rajchman (due to $\hat{\mu}(3n)=\hat{\mu}(n)$, $n\in\Z$). As in this last example, the obstructions for a measure to be Rajchman are often seen to be of arithmetical nature. The present work goes in this direction.  

\medskip
\noindent
As it concerns $t\rightarrow+\infty$, the Rajchman character of a measure $\mu$ on $\R$ is an information of local regularity. As is well-known, it says for example that $\mu$ has no atom; if $\hat{\mu}\in L^2(\R)$, then $\mu$ is absolutely continuous with respect to ${\cal L}_{\R}$ with an $L^2(\R)$ density; if $\hat{\mu}$ has some polynomial decay at infinity, one gets a lower bound on the Hausdorff dimension of $\mu$; etc. The Rajchman character can be reformulated as an equidistribution property modulo 1. Since $\hat{\mu}(t)\rightarrow0$ is equivalent to $\hat{\mu}(mt)\rightarrow0$ for any integer $m\not=0$, if $X$ is a real random variable with law $\mu$, then $\mu$ is Rajchman if and only if the law of $tX\mod 1$ converges, as $t\rightarrow+\infty$, to Lebesgue measure ${\cal L}_{\T}$ on $\T$.

\medskip
{\it Self-similar measures.} We now recall standard notions about self-similar measures on the real line $\R$, with a probabilistic point of view. We write ${\cal L}(X)$ for the law of a real random variable $X$. Let $N\geq0$ and real affine maps $\varphi_k(x)=r_kx+b_k$, with $r_k>0$, for $0\leq k\leq N$, and at least one $r_k<1$. We shall talk of ``strict contractions" in the case when $0<r_k<1$, for all $0\leq k\leq N$. This assumption will be considered principally in the second half of the article. For the sequel, we introduce the vectors $r=(r_k)_{0\leq k\leq N}$ and $b=(b_k)_{0\leq k\leq N}$. 

\smallskip
\noindent
Notice for what follows that for $n\geq0$, a composition  $\varphi_{k_{n-1}}\circ \cdots \circ\varphi_{k_0}$ has the explicit expression~:

$$\varphi_{k_{n-1}}\circ \cdots \circ\varphi_{k_0}(x)=r_{k_{n-1}}\cdots r_{ k_0}x+\sum_{l=0}^{n-1}b_{k_l}r_{k_{n-1}}\cdots r_{k_{l+1}}.$$

\medskip
\noindent
Consider the convex set ${\cal C}_N=\{p=(p_0,\cdots,p_N)~|~\forall j,~p_j>0,~\sum_jp_j=1\}$, open for the topology of the affine hyperplane $\{\sum_jp_j=1\}$. We denote its closure by $\bar{{\cal C}}_N$. Define :

$${\cal D}_N(r)=\{p\in\bar{{\cal C}}_N~|~\sum_{0\leq j\leq N}p_j\log r_j<0\}.$$

\smallskip
\noindent
This is a non-empty open subset of $\bar{{\cal C}}_N$, for the relative topology. Notice that ${\cal D}_N(r)=\bar{{\cal C}}_N$, in the case when the $(\varphi_k)_{0\leq k\leq N}$ are strict contractions.
 
\medskip
\noindent
Fixing a probability vector $p\in{\cal D}_N(r)$, we now compose the contractions at random, independently, according to $p$. Precisely, let $X_0$ be any real random variable and $(\varepsilon_n)_{n\geq0}$ be independent and identically distributed ($i.i.d.$) random variables, independent from $X_0$, and with law $p$, in other words $\P(\varepsilon_0=k)=p_k$, $0\leq k\leq N$. We consider the Markov chain $(X_n)_{n\geq0}$ on $\R$ defined by~:

$$X_{n}=\varphi_{\varepsilon_{n-1}}\circ \cdots \circ\varphi_{\varepsilon_{0}}(X_{0}),~n\geq0.$$

\medskip
\noindent
The condition $p\in{\cal D}_N(r)$, of contraction on average, can be rewritten as $\E(\log r_{\varepsilon_0})<0$. It implies that $(X_n)_{n\geq0}$ has a unique stationary (time invariant) measure, written as $\nu$. This follows for example from the fact that ${\cal L}(X_n)={\cal L}(Y_n)$, where~:

$$Y_n:=\varphi_{\varepsilon_{0}}\circ \cdots \circ\varphi_{\varepsilon_{n-1}}(X_0)=r_{\varepsilon_{0}}\cdots r_{ \varepsilon_{n-1}}X_0+\sum_{l=0}^{n-1}b_{\varepsilon_l}r_{\varepsilon_{0}}\cdots r_{\varepsilon_{l-1}}.$$

\noindent
As usual, $(Y_n)$ is more stable than $(X_n)$. Since $n^{-1}\log(r_{\varepsilon_0}\cdots r_{\varepsilon_{n-1}})\rightarrow\E(\log r_{\varepsilon_0})<0$, a.-s., as $n\rightarrow+\infty$, by the Law of Large Numbers, we get that $Y_n$ converges a.-s., as $n\rightarrow+\infty$, to~:

$$X:=\sum_{l\geq0}b_{\varepsilon_l}r_{\varepsilon_0}\cdots r_{\varepsilon_{l-1}}.$$

\smallskip
\noindent
Setting $\nu={\cal L}(X)$, we obtain that ${\cal L}(X_n)$ weakly converges to $\nu$. By construction, we have ${\cal L}(X_{n+1})=\sum_{0\leq j\leq N}p_j{\cal L}(X_n)\circ \varphi_j^{-1}$. Taking the limit as $n\rightarrow+\infty$, the measure $\nu$ verifies~:

\begin{equation}
\label{nun}
\nu=\sum_{0\leq j\leq N}p_j\nu\circ \varphi_j^{-1}.
\end{equation}

\smallskip
The previous convergence implies that the solution of this ``stable fixed point equation" is unique among Borel probability measures. Also, $\nu$ has to be of pure type, i.e. either purely atomic or absolutely continuous with respect to ${\cal L}_{\R}$ or else singular continuous, since each term in its Radon-Nikodym decomposition with respect to ${\cal L}_{\R}$ verifies \eqref{nun}. A few remarks are in order~:

\medskip
\noindent
$i)$ If $p\in{\cal C}_N$, the measure $\nu$ is purely atomic if and only if the $\varphi_j$ have a common fixed point $c$, in which case $\nu$ is the Dirac mass at $c$. Indeed, consider the necessity and suppose that $\nu$ has an atom. Let $a>0$ be the maximal mass of an atom and $E$ the finite set of points having mass $a$. Fixing any $c\in E$, the relation $\nu(\{c\})=\sum_jp_j\nu(\{\varphi_j^{-1}(c)\})$ furnishes $\varphi_j^{-1}(c)\in E$, $0\leq j\leq N$. Hence $\varphi_j^{-n}(c)\in E$, $n\geq0$, for all $j$. If $\varphi_j\not=id$, then $\varphi_j^{-1}(c)=c$, the set $\{\varphi_j^{-n}(c),~n\geq0\}$ being infinite otherwise. If $\varphi_j=id$, it fixes all points.

\medskip
\noindent
$ii)$ The equation for a hypothetical density $f$ of $\nu$ with respect to ${\cal L}_{\R}$, coming from \eqref{nun}, is~:

$$f=\sum_{0\leq j\leq N}p_jr_j^{-1}f\circ \varphi_j^{-1}.$$

\smallskip
\noindent
This ``unstable fixed point equation" is difficult to solve directly. It is equivalently reformulated into the fact that $((r_{\varepsilon_{n-1}}^{-1}\cdots r_{\varepsilon_0}^{-1})f\circ \varphi^{-1}_{\varepsilon_{n-1}}\cdots\circ\varphi_{\varepsilon_0}^{-1}(x))_{n\geq0}$ is a non-negative martingale (for its natural filtration), for Lebesgue a.-e. $x\in\R$. Notice that when $f$ exists and is bounded, then $p_j\leq r_j$ for all $j$, because $p_jr_j^{-1}\|f\|_{\infty}=\|p_jr_j^{-1}f\circ \varphi_j^{-1}\|_{\infty}\leq \|f\|_{\infty}$ and $\|f\|_{\infty}\not=0$.

\medskip
\noindent
$iii)$ Let $f(x)=ax+b$ be an affine map, with $a\not=0$. With the same $p\in{\cal D}_N(r)$, consider the conjugate system $(\psi_j)_{0\leq j\leq N}$, with $\psi_j(x)=f\circ \varphi_j\circ f^{-1}(x)=r_jx+b(1-r_j)+ab_j$. It has an invariant measure $w={\cal L}(aX+b)$ verifying the relation $\hat{w}(t)=\hat{\nu}(at)e^{2i\pi tb}$, $t\in\R$. In particular, $\nu$ is Rajchman if and only if $w$ is Rajchman.

\medskip
\noindent
$iv)$ When supposing that the $(\varphi_k)_{0\leq k\leq N}$ are strict contractions, some self-similar set $F$ can be introduced, where $F\subset\R$ is the unique non-empty compact set verifying the self-similarity relation $F=\cup_{0\leq k\leq N}\varphi_k(F)$. See for example Huchinson \cite{hut} for general properties of such sets. Introducing $\N=\{0,1,\cdots\}$ and the compact $S=\{0,\cdots,N\}^{\N}$, the hypothesis that the $(\varphi_k)_{0\leq k\leq N}$ are strict contractions implies that $F$ is a continuous (and even h\"olderian) image of $S$, in other words we have the following description~:

\begin{eqnarray}
F&=&\left\{{\sum_{l\geq0}b_{x_l}r_{x_0}\cdots r_{x_{l-1}},~(x_0,x_1,\cdots)\in S}\right\}.\nonumber
\end{eqnarray}

\medskip
\noindent 
Whereas in the general case a self-similar invariant measure can have $\R$ as topological support, when the $(\varphi_k)_{0\leq k\leq N}$ are strict contractions the compact self-similar set $F$ exists and supports any self-similar measure.

\bigskip
{\it Background and content of the article.} Coming back to the general case, we assume in the sequel that the $(\varphi_j)_{0\leq j\leq N}$ do not have a common fixed point (in particular $N\geq1$), so that $\mu$ is a continuous measure. A difficult problem is to characterize the absolute continuity of $\nu$ with respect to ${\cal L}_{\R}$ in terms of the parameters $r$, $b$ and $p$. An example with a long and well-known history is that of Bernoulli convolutions, corresponding to $N=1$, the affine contractions $\varphi_0(x)=\lambda x-1$, $\varphi_1(x)=\lambda x+1$, $0<\lambda<1$, and the probability vector $p=(1/2,1/2)$. Notice that when the $r_i$ are all equal (to some real in $(0,1)$), the situation is a little simplified, as $\nu$ is an infinite convolution (this is not true in general). Although we discuss below some works in this context, we will not present here the vast subject of Bernoulli convolutions, addressing the reader to detailed surveys, such as Peres-Schlag-Solomyak \cite{sixty} or Solomyak \cite{solo}.

\medskip
For general self-similar measures, an important aspect of the problem, that we shall not enter, and an active line of research, concerns the Hausdorff dimension of the measure $\nu$; cf the recent fundamental work of Hochman \cite{hochman} for example. In a large generality, see Falconer \cite{falco} and more recently Jaroszewska and Rams \cite{rams}, there is an $``\mbox{entropy}/\mbox{Lyapunov exponent}"$ upper-bound~:

$$\mbox{Dim}_{{\cal H}}(\nu)\leq\min\{1,s(p,r)\}\mbox{, where }s(p,r):=\frac{-\sum_{i=0}^Np_i\log p_i}{-\sum_{i=0}^Np_i\log r_i}.$$

\smallskip
\noindent
The quantity $s(p,r)$ is called the singularity dimension of the measure and can be $>1$. The equality $\mbox{Dim}_{{\cal H}}(\nu)=1$ does not mean that $\nu$ is absolutely continuous, but the inequality $s(p,r)<1$ surely implies that $\nu$ is singular. The interesting domain of parameters for the question of the absolute continuity of the invariant measure therefore corresponds to $s(p,r)\geq1$.

\medskip
We focus in this work on another fundamental tool, the Fourier transform $\hat{\nu}$. If $\nu$ is not Rajchman, the Riemann-Lebesgue lemma implies that $\nu$ is singular. This property was used by Erd\"os \cite{erdos} in the context of Bernoulli convolutions. He proved that if $1/2<\lambda<1$ is such that $1/\lambda$ is a Pisot number, then $\nu$ is not Rajchman. The reciprocal statement (for $1/2<\lambda<1$) was next shown by Salem \cite{salem}. As a result, for Bernoulli convolutions the Rajchman property always holds, except for a very particular set of parameters. Some works have next focused on the decay on average of the Fourier transform for general self-similar measures associated to strict contractions; cf Strichartz \cite{stri1,stri2}, Tsuji \cite{tsu}. In the same context, the non-Rajchman character was recently shown to hold for only a very small set of parameters by Li and Sahlsten \cite{lisa}, who showed that $\nu$ is Rajchman when $\log r_i/\log r_j$ is irrational for some $(i,j)$, with moreover some logarithmic decay of $\hat{\nu}$ at infinity, under a Diophantine condition. Next, Solomyak \cite{solo2} proved that outside a set of $r$ of zero Hausdorff dimension, $\hat{\nu}$ even has a power decay at infinity. 

\medskip
The aim of the present article is to study for general self-similar measures the exceptional set of parameters where the Rajchman property is not true, trying to follow the line of \cite{erdos} and \cite{salem}. We essentially show that $r$ and $b$ have to be closely related to some fixed Pisot number, as for Bernoulli convolutions. We first prove a general extension of the result of Salem \cite{salem}, reducing to a small island the set of parameters where the Rajchman property may not hold. Focusing then on this island of parameters, we provide a general characterization of the Rajchman character, appearing in this particular case as equivalent to absolute continuity with respect to ${\cal L}_{\R}$. Next, supposing that the $(\varphi_k)_{0\leq k\leq N}$ are strict contractions, we prove a partial extension of the theorem of Erd\"os \cite{erdos}, showing that for most parameters in the small island the Rajchman property is not true, with in general a few exceptions. We finally give some complements, first rather surprising numerical simulations involving the Plastic number, then an application to sets of uniqueness for trigonometric series.

\section{Statement of the results}

Let us place in the general context considered in the Introduction. Pisot numbers will play a central role in the analysis. Let us introduce a few definitions concerning Algebraic Number Theory; cf for example Samuel \cite{samuel} for more details.

\medskip

\begin{defi}

$ $

\noindent
A Pisot number is a real algebraic integer $\theta>1$, with conjugates (the other roots of its minimal unitary polynomial) of modulus strictly less than 1. Fixing such a $\theta>1$, denote its minimal polynomial as $Q=X^{s+1}+a_sX^s+\cdots+a_0\in\Z[X]$, of degree $s+1$, with $s\geq0$. If $s=0$, then $\theta$ is an integer $\geq2$. The images of $\mu\in\Q[\theta]$ by the $s+1$ $\Q$-homomorphisms $\Q[\theta]\rightarrow\C$ are the conjugates of $\mu$ corresponding to the field $\Q[\theta]$, in general denoted by $\mu=\mu^{(0)},\mu^{(1)},\cdots,\mu^{(s)}$.

\medskip
\noindent
$i)$ For $\alpha\in\Q[\theta]$, the trace $Tr_{\theta}(\alpha)$ is the trace of the linear operator $x\longmapsto \alpha x$ of multiplication by $\alpha$, considered from $\Q[\theta]$ to itself. As a general fact, $Tr_{\theta}(\alpha)\in\Q$.

\medskip
\noindent
$ii)$ Let $\Z[\theta]=\Z\theta^0+\cdots+\Z\theta^s$ be the subring generated by $\theta$ of the ring of algebraic integers of $\Q[\theta]$. We write ${\cal D}(\theta)$ for its $\Z$-dual (as a $\Z$-lattice), i.e.~:

$${\cal D}(\theta)=\{\alpha\in \Q[\theta],~Tr_{\theta}(\theta^n\alpha)\in\Z\mbox{, for }0\leq n\leq s\}.$$

\medskip
\noindent
It can be shown that ${\cal D}(\theta)=(1/Q'(\theta))\Z[\theta]$. As a classical fact, $Tr_{\theta}(\theta^n\alpha)\in\Z$, for all $n\geq0$, if this holds for $0\leq n\leq s$. Define~:

$${\cal T}(\theta)=\{\alpha\in\Q[\theta],~Tr_{\theta}(\theta^n\alpha)\in\Z\mbox{, for large }n\geq0\}.$$

\smallskip
\noindent
Then $\displaystyle {\cal T}(\theta)=\cup_{n\geq0}\theta^{-n}{\cal D}(\theta)=\frac{\Z[\theta,1/\theta]}{Q'(\theta)}$, with $\Z[\theta,1/\theta]$ the subring of $\Q[\theta]$ generated by $\{\theta,1/\theta\}$. 
\end{defi}

\medskip
\noindent
\begin{remark} In the context of the previous definition, introduce the integer-valued $(s+1)\times(s+1)$-companion matrix $M$ of $Q$~:

$$M=\({\begin{array}{cccc}0&1&\cdots&0\\
\vdots&\ddots&\ddots&\vdots\\
\vdots&\vdots&0&1\\
-a_0&\cdots&-a_{s-1}&-a_s\end{array}}\).$$

\medskip
\noindent
One may show that for any $\mu\in\Q[\theta]$, setting $V=(Tr_{\theta}(\theta^0\mu),\cdots,Tr_{\theta}(\theta^s\mu))$, then $\mu\in {\cal T}(\theta)$ if and only if there exists $n\geq0$ such that $VM^n$ has integral entries. 
\end{remark}

\medskip
We introduce special families of affine maps, that will play the role of canonical models for the analysis of the Rajchman property.

\begin{defi}

$ $

\noindent
Let $N\geq1$. A family of real affine maps $(\varphi_k)_{0\leq k\leq N}$ is said in Pisot form, if there exist a Pisot number $1/\lambda>1$, relatively prime integers $(n_k)_{0\leq k\leq N}$ and $\mu_k\in {\cal T}(1/\lambda)$, $0\leq k\leq N$, such that $\varphi_j(x)=\lambda^{n_j}x+\mu_j\mbox{, for all }0\leq j\leq N$.

\end{defi}

\noindent
\begin{remark} 
If $(\varphi_j)_{0\leq j\leq N}$ is in Pisot form, then the $(\lambda,(n_j),(\mu_j))$ are uniquely determined. Indeed, if the $(\lambda',(n'_j),(\mu'_j))$ also convene, we just need to show that $\lambda=\lambda'$. Taking some collection of integers $(a_j)$ realizing a Bezout relation $1=\sum_j a_jn_j$, we have~:

$$\lambda=\lambda^{\sum_ja_jn_j}=\lambda'^{\sum_ja_jn'_j}=\lambda'^p,$$

\noindent
for some $p\geq1$. Idem, $\lambda'=\lambda^q$, for some $q\geq1$. Hence $pq=1$, giving $p=q=1$ and $\lambda=\lambda'$. 

\end{remark}

\medskip
As a first result, extending \cite{salem}, the analysis of the non-Rajchman character of the invariant measure requires to consider families in Pisot form. 

\begin{thm}
\label{part1}

$ $

\noindent
Let $N\geq 1$, $p\in{\cal C}_N$ and affine maps $\varphi_k(x)=r_kx+b_k$, $r_k>0$, for $0\leq k\leq N$, with no common fixed point, and $\sum_{0\leq j\leq N}p_j\log r_j<0$. The invariant measure $\nu$ is not Rajchman if and only if there exists $f(x)=ax+b$, $a\not=0$, such that the conjugate system $(f\circ \varphi_j\circ f^{-1})_{0\leq j\leq N}$ is in Pisot form, for some Pisot number $1/\lambda>1$, with invariant measure $w$ verifying $\hat{w}(\lambda^{-k})\not\rightarrow_{k+\infty}0$.

\end{thm}

\noindent
In particular, one gets that $r_j=\lambda^{n_j}$, for all $j$, for some Pisot number $1/\lambda>1$ et relatively prime integers $(n_k)_{0\leq k\leq N}$. Hence, up to an affine change of variables, the non-Rajchman character of the invariant measure $\nu$ can be read on the sequence $(\lambda^{-k})_{k\geq0}$, as in \cite{erdos}. In a second step, we provide a general analysis of families in Pisot form.

\medskip
We now fix a Pisot number $1/\lambda>1$, an integer $N\geq 1$, relatively prime integers $(n_k)_{0\leq k\leq N}$ and $(\mu_k)_{0\leq k\leq N}\in ({\cal T}(1/\lambda))^{N+1}$, such that $\varphi_k(x)=\lambda^{n_k}x+\mu_k$, for $0\leq k\leq N$. Let $p\in{\cal C}_N$ be such that $\sum_{0\leq j\leq N}p_jn_j>0$ and $i.i.d.$ random variables $(\varepsilon_n)_{n\in\Z}$, with $\P(\varepsilon_0=k)=p_k$, $0\leq k\leq N$. We introduce cocycle notations $(S_l)_{l\in\Z}$, where $S_0=0$ and for $l\geq1$~:

$$S_l=n_{\varepsilon_0}+\cdots+n_{\varepsilon_{l-1}},~~~S_{-l}=-n_{\varepsilon_{-l}}-\cdots-n_{\varepsilon_{-1}}.$$

\medskip
\noindent
An important preliminary remark is that when $\mu\in{\cal T}(1/\lambda)$ and $k\geq0$ is large enough, we have~:

$$\lambda^{-k}\mu+\sum_{1\leq j\leq s}\alpha_j^k\mu^{(j)}=Tr_{1/\lambda}(\lambda^{-k}\mu)\in\Z,$$

\noindent
where the $(\alpha_j)_{0\leq j\leq s}$ are the conjugates of $1/\lambda=:\alpha_0$ and the $(\mu^{(j)})_{0\leq j\leq s}$ that of $\mu=\mu^{(0)}$, corresponding to the field $\Q[\lambda]$. Since $|\alpha_j|<1$, for $1\leq j\leq s$, and $(S_l)$ is a.-s. transient with a non-zero linear speed to $-\infty$, as $l\rightarrow-\infty$, by the Law of Large Numbers, this ensures that for any $k\in\Z$, the random variable $\sum_{l\in\Z}\mu_{\varepsilon_l}\lambda^{k+S_l}\mod 1$ is well-defined as a $\T$-valued random variable.

\medskip
In the sequel we use standard inner products and Euclidean norms on all spaces $\R^n$.

\begin{thm}
\label{part2}

$ $

\noindent
Let $1/\lambda>1$ be a Pisot number of degree $s+1$. Let $N\geq 1$, relatively prime integers $(n_k)_{0\leq k\leq N}$ and $(\mu_k)_{0\leq k\leq N}\in ({\cal T}(1/\lambda))^{N+1}$, such that $\varphi_k(x)=\lambda^{n_k}x+\mu_k$, $0\leq k\leq N$. Let $p\in{\cal C}_N$ be such that $\sum_{0\leq j\leq N}p_jn_j>0$ and $i.i.d.$ random variables $(\varepsilon_n)_{n\in\Z}$, with law $p$. Let $(S_l)_{l\in\Z}$ be the cocycle notations associated to the $(n_{\varepsilon_i})_{i\in\Z}$. The real random variable $X=\sum_{l\geq0}\mu_{\varepsilon_l}\lambda^{S_l}$ has law $\nu$.

\medskip
\noindent
i) Let the $\T$-valued random variables $Z_k=\sum_{l\in\Z}\mu_{\varepsilon_l}\lambda^{k+S_l}$, $k\in\Z$. Then $\lambda^{-n}X\mod 1$ converges, as $n\rightarrow+\infty$, to a probability measure $m$ on $\T$, verifying, for all $f\in C(\T,\R)$ and all $k\in\Z$~:

$$\int_{\T} f(x)~dm(x)=\frac{1}{\E(n_{\varepsilon_0})}\sum_{0\leq r<n^*}\E\left[{f\({Z_{k+r}}\)1_{S_{-u}<-r,u\geq1}}\right],$$ 

\noindent
where $n^*=\max_{0\leq k\leq N}n_k$. More generally, $\lambda^{-n}(X,\lambda^{-1}X,\cdots,\lambda^{-s}X) \mod\Z^{s+1}$ converges in law, as $n\rightarrow+\infty$, to a probability measure $\cal{M}$ on $\T^{s+1}$, with one-dimensional marginals $m$, verifying~: 

$$\int_{\T^{s+1}} f(x)~d{\cal M}(x)=\frac{1}{\E(n_{\varepsilon_0})}\sum_{0\leq r<n^*}\E\left[{f\({Z_{k+r},Z_{k+r-1},\cdots,Z_{k+r-s}}\)1_{S_{-u}<-r,u\geq1}}\right],$$ 

\smallskip
\noindent
for all $f\in C(\T^{s+1},\R)$ and all $k\in\Z$.

\medskip
\noindent
ii) If the $(\varphi_k)_{0\leq k\leq N}$ do not have a common fixed point (i.e. if $\nu$ is continuous), denoting by $Z$ a $\T^{s+1}$-valued random variable with law $\cal{M}$, then for any $0\not=n=(n_0,\cdots,n_s)^t\in\Z^{s+1}$, $\langle Z,n\rangle$ has a continuous law; in particular, $m$ and ${\cal M}$ are continuous measures. If the $(\varphi_k)_{0\leq k\leq N}$ have a common fixed point, there exists a rational number $p/q$ such that $m=\delta_{p/q}$ and ${\cal M}=(\delta_{p/q})^{\otimes (s+1)}$.

\medskip
\noindent
iii) Either ${\cal M}\perp{\cal L}_{\T^{s+1}}$ or ${\cal M}={\cal L}_{\T^{s+1}}$. Also, ${\cal M}={\cal L}_{\T^{s+1}}\Leftrightarrow \nu\mbox{ is Rajchman }\Leftrightarrow\nu\ll{\cal L}_{\R}$.

 \end{thm}
 
In the context of the previous theorem, $\nu$ and ${\cal M}$ are always of the same nature, with respect to the uniform measure of the space they live on. In particular, ${\cal M}$ is also of pure type. We finally consider families in Pisot form, when supposing that the $(\varphi_k)_{0\leq k\leq N}$ are strict contractions.

 \begin{thm}
 \label{part3}

$ $
 
 \noindent
Let $N\geq 1$ and $\varphi_k(x)=\lambda^{n_k}x+\mu_k$, for $0\leq k\leq N$, with $1/\lambda>1$ a Pisot number, relatively prime integers $(n_k)_{0\leq k\leq N}$, with $n_k\geq1$ and $\mu_k\in {\cal T}(1/\lambda)$, for $0\leq k\leq N$. When $p\in{\cal C}_N$ is fixed, we denote by $m$ the measure on $\T$ of Theorem \ref{part2}, $i)$.

\medskip
\noindent
i) For any $p\in{\cal C}_N$, if the invariant measure $\nu$ is Rajchman, then it is absolutely continuous with respect to ${\cal L}_{\R}$, with a density bounded and with compact support.

\medskip
\noindent
ii) There exists $0\not=a\in\Z$ such that for any $k\not=0$, for any $p\in{\cal C}_N$ outside finitely many real-analytic graphs of dimension $\leq N-1$ (points if $N=1$), we have $\hat{m}(ak)\not=0$. In this case, $m\not={\cal L}_{\T}$ and $\nu$ is not Rajchman.

\end{thm}

\medskip
\noindent
\begin{remark}
In Theorem \ref{part3} $ii)$, observe that when making $k$ vary, we obtain that for all $p\in{\cal C}_N$ outside a countable number of real-analytic graphs of dimension less than or equal to $ N-1$ (points if $N=1$), then $\hat{m}(ak)\not=0$, for all $k\in\Z$. Part $ii)$ of Theorem \ref{part3} relies on an indirect argument, based on the analysis of the regularity of $\hat{m}(n)$, for some fixed $n\in\Z$, as a function of $p\in{\cal C}_N$. 
\end{remark}

\medskip
\noindent
\begin{remark}
On the existence of singular measures in the non-homogeneous case, we are essentially aware of the non-explicit examples, using algebraic curves, of Neunh\"auserer \cite{neun}. As suggested by the referee, Theorem \ref{part3} allows to give in the non-homogeneous case an explicit example of a continuous singular and not Rajchman invariant measure $\nu$ with singularity dimension $>1$. Indeed, take for $1/\lambda>1$ the Plastic number, i.e. the real root of $X^3-X-1$. This is the smallest Pisot number; cf Siegel \cite{siegel}. We have $1/\lambda=1.3247...$. Let $N=1$ and $\varphi_0(x)=\lambda x,~\varphi_1(x)=\lambda^2x+1$. For $p=(p_0,p_1)\in{\cal C}_1$, if $\nu$ is absolutely continuous with respect to ${\cal L}_{\R}$, then, by Theorem \ref{part3} $i)$, the density has to be bounded. By remark $ii)$ in the Introduction, this implies that $p_0\leq \lambda=0.7548...$ and $p_1\leq\lambda^2$. Now, as detailed in the last section, the similarity dimension in this case is $>1$ if and only if $0,203...<p_0<0,907..$. For example we can conclude that for $p_0\in[0.76,0.90]$, the measure $\nu$ is continuous, singular with respect to ${\cal L}_{\R}$, not Rajchman and with similarity dimension $>1$. Still for the system $\varphi_0(x)=\lambda x$, $\varphi_1(x)=\lambda^2x+1$, we will give in the last section a strong numerical support for the fact that $\nu$ is in fact continuous singular and not Rajchman for all $p\in{\cal C}_1$. \end{remark}

\medskip
\noindent
\begin{remark} In the context of Theorem \ref{part3}, it would be important to determine all the exceptional parameters where $\nu$ is absolutely with respect to ${\cal L}_{\R}$. Let us give some examples where the exceptional set in Theorem \ref{part3} $ii)$ is non-empty~:

\medskip
\noindent
1) Let $1/\lambda=N\geq1$ and $\varphi_k(x)=(x+k)/(N+1)$, with $p_k=1/(N+1)$, for $0\leq k\leq N$; then $\nu$ is Lebesgue measure on $[0,1]$. 

\medskip
\noindent
2) Take for $1/\lambda>1$ the Plastic number, $N=1$ and this time $\varphi_0(x)=\lambda^2x,~\varphi_1(x)=\lambda^3x+1$. One may verify that the similarity dimension is $<1$ for all $p\in{\cal C}_1$, except for $p=(\lambda^2,\lambda^3)$, where it equals one. Thus the invariant measure $\nu$ is singular for $p\in{\cal C}_1$ with $p\not=(\lambda^2,\lambda^3)$. Another way, if $\nu$ is absolutely continuous with respect to ${\cal L}_{\R}$, then its density has to be bounded by Theorem \ref{part3}. Therefore, $p_0\leq\lambda^2$ and $p_1\leq \lambda^3$, using remark $ii)$ in the Introduction. Since $\lambda^2+\lambda^3=1$, we have $p_0=\lambda^2$ and $p_1=\lambda^3$. As a result, when $p=(p_0,p_1)\not=(\lambda^2,\lambda^3)$ and $p_0>0$, $p_1>0$, then $\nu$ is continuous singular and not Rajchman. When $p=(\lambda^2,\lambda^3)$, set $I=[0,1+\lambda]$ and notice that $\varphi_0(I)=[0,1]$, $\varphi_1(I)=[1,1+\lambda]$. Hence, Lebesgue a.-e.~:

$$1_{I}=1_{\varphi_0(I)}+1_{\varphi_1(I)}=p_0\lambda^{-2}1_{\varphi_0(I)}+p_1\lambda^{-3}1_{\varphi_1(I)},$$

\smallskip
\noindent
meaning that $\nu=\frac{1}{1+\lambda}{\cal L}_{I}$. Taking for $1/\lambda$ the supergolden ratio (the real root of $X^3-X^2-1$; the fourth Pisot number), one gets the same situation with the system $(\lambda x+1,\lambda^3x)$, the exceptional parameters being then $(\lambda,\lambda^3)$, giving for $\nu$ the uniform probability measure on $[0,\lambda^{-3}]$. 

\medskip
\noindent
3) When $1/\lambda>1$ is the Plastic number, $N=2$, $\varphi_0(x)=\lambda^2x$, $\varphi_1(x)=\lambda^3x+1$, $\varphi_2(x)=\lambda^3x+1$ and $p_0=\lambda^2$, $p_1=\lambda^3\alpha$, $p_2=\lambda^3(1-\alpha)$, then $\nu=\frac{1}{1+\lambda}{\cal L}_{[0,1+\lambda]}$, for all $0<\alpha<1$. This is an example, a little degenerated, of a one-dimensional real-analytic graph where the corresponding invariant measure $\nu$ is absolutely continuous with respect to ${\cal L}_{\R}$.

\medskip
It would be interesting to find more developed examples, where $\nu$ is absolutely continuous with respect to ${\cal L}_{\R}$. A difficulty is that a priori the probability vector $p$ has to be chosen in accordance with the polynomial equations verified by $\lambda$.  \end{remark}

\section{Proof of Theorem \ref{part1}}

Let $N\geq1$ and $(\varphi_k)_{0\leq k\leq N}$, with $\varphi_k(x)=r_kx+b_k$, where $r_k>0$, and having no common fixed point. Fixing $p\in {\cal C}_N$, introduce $i.i.d.$ random variables $(\varepsilon_n)_{n\geq0}$ with law $p$, to which $\P$ and $\E$ refer. By hypothesis, $\E(\log r_{\varepsilon_0})<0$. Recall that the invariant measure $\nu$ is the law of the random variable $\sum_{l\geq0}b_{\varepsilon_l}r_{\varepsilon_0}\cdots r_{\varepsilon_{l-1}}$ and that $\nu$ is supposed to be non Rajchman. Without loss of generality, we assume that $0<r_0\leq r_1\leq\cdots\leq r_N$, with therefore $r_0<1$.

\medskip
The proof has three parts. First we show that $\log r_i/\log r_j\in\Q$, for all $0\leq i\not=j\leq N$. From this, we will get that $r_j=\lambda^{n_j}$, for some $0<\lambda<1$ and integers $(n_j)$. We then show that the non Rajchman character of $\nu$ can be seen on a subsequence of the form $(\alpha\lambda^{-k})_{k\geq0}$. We finally prove that $1/\lambda$ is a Pisot number and the family $(\varphi_k)_{0\leq k\leq N}$ is affinely conjugated with one in Pisot form.

\medskip
\noindent
{\it Step 1.}  Let us show that if ever $\log r_i/\log r_j\not\in\Q$, for some $0\leq i\not=j\leq N$, then $\nu$ is Rajchman. This is established in \cite{lisa} for strict contractions. We simplify their proof. 

\medskip
 For $n\geq1$, consider the random walk $S_n=-\log r_{\varepsilon_0}-\cdots-\log r_{\varepsilon_{n-1}}$, with $S_0=0$. For a real $s\geq0$, introduce the finite stopping time $\tau_s=\min\{n\geq0,~S_n>s\}$ and write ${\cal T}_s$ for the corresponding sub-$\sigma$-algebra of the underlying $\sigma$-algebra. Taking $\alpha>0$ and $s\geq0$~:

\begin{eqnarray}
\hat{\nu}(\alpha e^s)&=&\E\({e^{2\pi i\alpha e^s\sum_{l\geq0}b_{\varepsilon_l}e^{-S_l}}}\)\nonumber\\
&=&\E\({e^{2\pi i\alpha e^{s}\sum_{0\leq l<\tau_s}b_{\varepsilon_l}e^{-S_l}}e^{2\pi i\alpha e^{-S_{\tau_s}+s}\sum_{l\geq\tau_s}b_{\varepsilon_l}e^{-S_l+S_{\tau_s}}}}\).\nonumber\end{eqnarray}

\smallskip
\noindent
In the expectation, the first exponential term is ${\cal T}_s$-measurable. Also, the conditional expectation of the second exponential term with respect to ${\cal T}_s$ is just $\hat{\nu}(\alpha e^{-S_{\tau_s}+s})$, as a consequence of the strong Markov property. It follows that~:

$$\hat{\nu}(\alpha e^s)=\E\({\hat{\nu}(\alpha e^{-S_{\tau_s}+s})e^{2\pi i\alpha e^{s}\sum_{0\leq l<\tau_s}b_{\varepsilon_l}e^{-S_l}}}\).$$

\medskip
\noindent
This gives $|\hat{\nu}(\alpha e^s)|\leq\E\({|\hat{\nu}(\alpha e^{-S_{\tau_s}+s})|}\)$, so by the Cauchy-Schwarz inequality and the Fubini theorem, which directly applies, consecutively~:

\begin{eqnarray}
|\hat{\nu}(\alpha e^s)|^2\leq\E\({|\hat{\nu}(\alpha e^{-S_{\tau_s}+s})|^2}\)&=&\E\({\int_{\R^2}e^{2\pi i\alpha e^{-S_{\tau_s}+s}(x-y)}~d\nu(x)d\nu(y)}\)\nonumber\\
&=&\int_{\R^2}\E\({e^{2\pi i\alpha e^{-S_{\tau_s}+s}(x-y)}}\)~d\nu(x)d\nu(y)\nonumber\\
&\leq&\int_{\R^2}\left|{\E\({e^{2\pi i\alpha e^{-S_{\tau_s}+s}(x-y)}}\)}\right|~d\nu(x)d\nu(y).\nonumber\end{eqnarray}

\medskip
\noindent
Let $Y:=-\log r_{\varepsilon_0}$. The law of $Y$ is non-lattice, since some $\log r_i/\log r_j\not\in\Q$ and $p_k>0$ for all $0\leq k\leq N$. As $Y$ is integrable, with $0<\E(Y)<\infty$, it is a well-known consequence of the Blackwell theorem on the law of the overshoot that (see for instance Woodroofe \cite{wood}, chap. 2, thm 2.3), that~:

$$\E(g(S_{\tau_s}-s))\rightarrow \frac{1}{\E(S_{\tau_0})}\int_0^{+\infty}g(x)\P(S_{\tau_0}>x)~dx\mbox{, as }s\rightarrow+\infty,$$

\noindent
for any Riemann-integrable $g$ defined on $\R^+$. Here, all $S_{\tau_s}-s$, for $s\geq0$, (and in particular $S_{\tau_0}$) have support in some $[0,A]$. Therefore, $\P(S_{\tau_0}>x)=0$ for large $x>0$. For any $\alpha>0$, by dominated convergence (letting $s\rightarrow+\infty$)~:

$$\limsup_{t\rightarrow+\infty}|\hat{\nu}(t)|^2\leq  \frac{1}{\E(S_{\tau_0})}\int_{\R^2}\left|{\int_{0}^{+\infty}e^{2\pi i\alpha e^{-u}(x-y)}\P(S_{\tau_0}>u)du}\right|d\nu(x)d\nu(y).$$

\medskip
\noindent
The inside term (in the modulus) is uniformly bounded with respect to $(x,y)\in\R^2$. We shall use dominated convergence once more, this time with $\alpha\rightarrow+\infty$. It is sufficient to show that for $\nu^{\otimes2}$-almost every $(x,y)$, the inside term goes to zero. Since the measure $\nu$ is non-atomic, $\nu^{\otimes2}$-almost-surely, $x\not=y$. If for example $x>y$~:

$$\int_{0}^{+\infty}e^{2\pi i\alpha e^{-u}(x-y)}\P(S_{\tau_0}>u)du=\int_0^{x-y}e^{2\pi i\alpha t}\P(S_{\tau_0}>\log((x-y)/t)~\frac{dt}{t},$$

\medskip
\noindent
making the change of variable $t=e^{-u}(x-y)$. The last integral now converges to 0, as $\alpha\rightarrow+\infty$, by the Riemann-Lebesgue lemma. Hence, $\lim_{t\rightarrow+\infty}\hat{\nu}(t)=0$. This ends the proof of this step.

\bigskip
\noindent
{\it Step 2.} As $\nu$ is not Rajchman, from {\it Step 1}, $\log r_i/\log r_j\in\Q$, for all $(i,j)$. Hence $r_j=r_0^{p_j/q_j}$, with integers $p_j\in\Z$, $q_j\geq1$, for $1\leq j\leq N$. Let :

$$n_0=\prod_{1\leq l\leq N}q_l\geq1\mbox{ and }n_j=p_j\prod_{1\leq l\leq N,l\not=j}q_l\in\Z,~1\leq j\leq N.$$

\noindent
Recall that $0<r_0<1$. Setting $\lambda=r_0^{1/n_0}\in(0,1)$, one has $r_j=\lambda^{n_j}$, $0\leq j\leq N$. Up to taking some positive integral power of $\lambda$, one can assume that $\gcd(n_0,\cdots,n_N)=1$. Recall in passing that the set of Pisot numbers is stable under positive integral powers. The condition $\E(\log r_{\varepsilon_0})<0$ rewrites into $\E(n_{\varepsilon_0})>0$ and we have $n_N\leq\cdots\leq n_0$, with $n_0\geq1$.

\medskip
Using some sub-harmonicity, we shall now show that one can reinforce the assumption that $\hat{\nu}(t)$ is not converging to $0$, as $t\rightarrow+\infty$.

\begin{lemme}
\label{subhar}

$ $

\noindent
There exists $1\leq\alpha\leq1/\lambda$ and $c>0$ such that $\hat{\nu}(\alpha\lambda^{-k})=c_ke^{2i\pi\theta_k}$, written in polar form, verifies $c_k\rightarrow c$, as $k\rightarrow+\infty$. 
\end{lemme}

\noindent
{\it Proof of the lemma~:}

\noindent
Let us write this time $S_n=n_{\varepsilon_0}+\cdots+n_{\varepsilon_{n-1}}$, for $n\geq1$, with $S_0=0$. Since $\E(n_{\varepsilon_0})>0$, $(S_n)$ is transient to $+\infty$. Introduce the random ladder epochs $0=\sigma_0<\sigma_1<\cdots$, where inductively $\sigma_{k+1}$ is the first time $n\geq0$ with $S_n>S_{\sigma_{k}}$. Let $S'_k=S_{\sigma_{k}}$. The $(S'_k-S'_{k-1})_{k\geq1}$ are $i.i.d.$ random variables with law ${\cal L}(S_{\tau_0})$ and support in $\{1,\cdots,n_0\}$. Since $\gcd(n_0,\cdots,n_N)=1$, the support of the law of $S_{\tau_0}$ generates $\Z$ as an additive group (cf for example Woodroofe \cite{wood}, thm 2.3, second part). For an integer $u\geq1$ large enough, we can fix integers $r\geq1$ and $s\geq1$ such that the support of the law of $S'_r$ contains $u$ and that of $S'_s$ contains $u+1$, both supports being included in some $\{1,\cdots,M\}$, with therefore $1\leq u\leq u+1\leq M$. Proceeding as in {\it Step 1}, for any $t\in\R$~:

$$\hat{\nu}(t)=\E\({e^{2\pi it\sum_{l\geq0}b_{\varepsilon_l}\lambda^{S_l}}}\)=\E\({\hat{\nu}(t\lambda^{S'_{r}})e^{2\pi it\sum_{0\leq l<\sigma_r}b_{\varepsilon_l}\lambda^{S_l}}}\).$$

\medskip
\noindent
Doing the same thing with $S'_s$ and taking modulus gives~:

\begin{equation}
\label{itera}
|\hat{\nu}(t)|\leq\E\({|\hat{\nu}(t\lambda^{S'_{r}})|}\)\mbox{ and }|\hat{\nu}(t)|\leq\E\({|\hat{\nu}(t\lambda^{S'_{s}})|}\).
\end{equation}
 
\smallskip
\noindent 
In particular, $|\hat{\nu}(t)|\leq\max_{1\leq l\leq M}|\hat{\nu}(\lambda^lt)|$. We now set~:

$$V_{\alpha}(k):=\max_{k\leq l<k+M}|\hat{\nu}(\alpha\lambda^l)|,~k\in\Z,~\alpha>0.$$

\noindent
The previous remarks imply that $V_{\alpha}(k)\leq V_{\alpha}(k+1)$, $k\in\Z$, $\alpha>0$.

\medskip
Since $\nu$ is not Rajchman, $|\hat{\nu}(t_l)|\geq c'>0$, along some sequence $t_l\rightarrow+\infty$. Write $t_l=\alpha_l\lambda^{-k_l}$, with $1\leq \alpha_l\leq 1/\lambda$ and $k_l\rightarrow+\infty$. Up to taking a subsequence, $\alpha_l\rightarrow\alpha\in[1,1/\lambda]$. Fixing $k\in\Z$~:

$$c'\leq V_{\alpha_l}(-k_l)\leq V_{\alpha_l}(-k),$$

\smallskip
\noindent
as soon as $l$ is large enough. By continuity, letting $l\rightarrow+\infty$, we get $c'\leq V_{\alpha}(-k)$, $k\in\Z$. As $k\longmapsto V_{\alpha}(-k)$ is non-increasing, $V_{\alpha}(-k)\rightarrow c\geq c'$, as $k\rightarrow+\infty$. We now show that necessarily $|\hat{\nu}(\alpha\lambda^{-k})|\rightarrow c$, as $k\rightarrow+\infty$.

\medskip
If this were not true, there would exist $\varepsilon>0$ and $(m_k)\rightarrow+\infty$, with $|\hat{\nu}(\alpha\lambda^{-m_k})|\leq c-\varepsilon$. Using $V_{\alpha}(-k)\rightarrow c$ and $|\hat{\nu}(\alpha\lambda^{-m_k})|\leq c-\varepsilon$, as $k\rightarrow+\infty$, consider \eqref{itera} with $r$ and $t=\alpha\lambda^{-m_k-u}$ and next with $s$ and $t=\alpha\lambda^{-m_k-u-1}$. Since $u$ is in the support of the law of $S'_r$ and $u+1$ is in the support of the law of $S'_s$, we obtain the existence of some $c_1<c$ such that for $k$ large enough~:

$$\max\{|\hat{\nu}(\alpha\lambda^{-m_k-u})|,|\hat{\nu}(\alpha\lambda^{-m_k-u-1})|\}\leq c_1<c.$$

\smallskip
\noindent
Again via \eqref{itera}, with successively $r$ and $t=\alpha\lambda^{-m_k-2u}$, next $r$ and $t=\alpha\lambda^{-m_k-2u-1}$ and finally $s$ and $t=\alpha\lambda^{-m_k-2u-2}$, still using that $u$ is in the support of the law of $S'_r$ and $u+1$ in the support of the law of $S'_s$, we get some $c_2<c$ such that for $k$ large enough~:

$$\max\{|\hat{\nu}(\alpha\lambda^{-m_k-2u})|,~|\hat{\nu}(\alpha\lambda^{-m_k-2u-1})|,~|\hat{\nu}(\alpha\lambda^{-m_k-2u-2})|\}\leq c_2<c.$$

\medskip\noindent
Etc, for some $c_{M-1}<c$ and $k$ large enough~:

$$\max\{|\hat{\nu}(\alpha\lambda^{-m_k-(M-1)u})|,~\cdots,~|\hat{\nu}(\alpha\lambda^{-m_k-(M-1)u-(M-1)})|\}\leq c_{M-1}<c.$$

\medskip\noindent
This contradicts the fact that $V_{\alpha}(-k)\rightarrow c$, as $k\rightarrow\infty$. We conclude that $|\hat{\nu}(\alpha\lambda^{-k})|\rightarrow c$, as $k\rightarrow\infty$, and this ends the proof of the lemma.

\fin

\bigskip
\noindent
{\it Step 3.} We complete the proof of Theorem \ref{part1}. In this part, introduce the notation $\|x\|=\mbox{dist}(x,\Z)$, for $x\in\R$. Let us consider any $1\leq\alpha\leq 1/\lambda$, with $\hat{\nu}(\alpha\lambda^{-k})=c_ke^{2i\pi\theta_k}$, verifying $c_k\rightarrow c>0$, as $k\rightarrow+\infty$. The existence of such a $\alpha$ was shown in {\it Step 2}. We start from the relation~:

$$\hat{\nu}(\alpha\lambda^{-k})=\sum_{0\leq j\leq N}p_je^{2i\pi\alpha\lambda^{-k}b_j}\hat{\nu}(\alpha\lambda^{-k+n_j}),$$

\noindent
obtained when conditioning with respect to the value of $\varepsilon_0$. This furnishes for $k\geq0$~:

$$c_k=\sum_{0\leq j\leq N}p_je^{2i\pi(\alpha\lambda^{-k}b_j+\theta_{k-n_j}-\theta_k)}c_{k-n_j}.$$

\noindent
We rewrite this as~:

$$\sum_{0\leq j\leq N}p_j\[{e^{2i\pi(\alpha\lambda^{-k}b_j+\theta_{k-n_j}-\theta_k)}-1}\]c_{k-n_j}=c_k-\sum_{0\leq j\leq N}p_jc_{k-n_j}=\sum_{0\leq j\leq N}p_j(c_k-c_{k-n_j}).$$

\noindent
Let $K>0$ be such that $c_{k-n_j}\geq c/2>0$, for $k\geq K$ and all $0\leq j\leq N$. For $L>n^*$, where $n^*=\max_{0\leq j\leq N}|n_j|$, we sum the previous equality on $K\leq k\leq K+L$~:

$$\sum_{0\leq j\leq N}p_j\sum_{k=K}^{K+L}c_{k-n_j}\[{e^{2i\pi(\alpha\lambda^{-k}b_j+\theta_{k-n_j}-\theta_k)}-1}\]=\sum_{0\leq j\leq N}p_j\({\sum_{k=K}^{K+L}c_k-\sum_{k=K-n_j}^{K+L-n_j}c_{k}}\).$$

\noindent 
Observe that the right-hand side involves a telescopic sum and is bounded by $2n^*$ (using that $|c_k|\leq 1$), uniformly in $K$ and $L$. In the left hand-hand side, we take the real part and use that $1-\cos(2\pi x)=2(\sin \pi x)^2$, which, as is well-known, has the same order as $\|x\|^2$. We obtain, for some constant $C$, that for $K$ and $L$ large enough~:

$$\frac{c}{2}\sum_{0\leq j\leq N}p_j\sum_{k=K}^{K+L}\|\alpha\lambda^{-k}b_j+\theta_{k-n_j}-\theta_k\|^2\leq C.$$

\noindent
Introducing the constants $p_*=\min_{0\leq j\leq N}p_j>0$ and $C'=2C/(cp_*)$, we get that for all $0\leq j\leq N$ and $K,L$ large enough~:

\begin{equation}
\label{pisot}
\sum_{k=K}^{K+L}\|\alpha\lambda^{-k}b_j+\theta_{k-n_j}-\theta_k\|^2\leq C'.
\end{equation}

\noindent
In the sequel, we distinguish two cases~: there exists a non-zero translation among the $(\varphi_k)_{0\leq k\leq N}$ (case 1) or not (case 2). 

\medskip
\noindent
- {\it Case 1.} For any non-zero-translation $\varphi_j(x)=x+b_j$, we have $n_j=0$ and $b_j\not=0$. Then \eqref{pisot} gives that for $K,L$ large enough~:

$$\sum_{k=K}^{K+L}\|\alpha\lambda^{-k}b_j\|^2\leq C'.$$

\smallskip
\noindent
This implies that $(\|\alpha b_j\lambda^{-k}\|)_{k\geq0}\in l^2(\N)$. By a classical theorem of Pisot, cf Cassels \cite{cassels}, chap. 8, Theorems I and II, we obtain that $1/\lambda$ is a Pisot number and $b_j=(1/\alpha)\mu_j$, with $\mu_j\in{\cal T}(1/\lambda)$. Consider now the non-translations $\varphi_j(x)=\lambda^{n_j}x+b_j$, $n_j\not=0$. By \eqref{pisot}, for any $r\geq0$ and $K,L$ large enough (depending on $r$)~:

$$\sum_{k=K}^{K+L}\|\alpha\lambda^{-k+rn_j}b_j+\theta_{k-(r+1)n_j}-\theta_{k-rn_j}\|^2\leq C'.$$

 \noindent
Fixing $l_j\geq1$ and summing over $0\leq r\leq l_j-1$, making use of the triangular inequality and of $(x_1+\cdots+x_n)^2\leq n(x_1^2+\cdots+x_n^2)$, we obtain, for $K,L$ large enough (depending on $l_j$)~:

\begin{equation}
\label{ver1}
\sum_{k=K}^{K+L}\left\|{\alpha\lambda^{-k}b_j\({\frac{1-\lambda^{l_jn_j}}{1-\lambda^{n_j}}}\)+\theta_{k-l_jn_j}-\theta_{k}}\right\|^2\leq l_jC'.
\end{equation}

\smallskip
\noindent
Changing $k$ into $k+l_jn_j$, we obtain, for $K,L$ large enough (depending on $l_j$)~:

\begin{equation}
\label{ver2}
\sum_{k=K}^{K+L}\left\|{\alpha\lambda^{-k}b_j\({\frac{1-\lambda^{-l_jn_j}}{1-\lambda^{n_j}}}\)+\theta_{k+l_jn_j}-\theta_{k}}\right\|^2\leq l_jC'.
\end{equation}

\smallskip
\noindent
Let $1=\sum_{0\leq j\leq N}l_jn_j$ be a Bezout relation and $J\subset\{0,\cdots,N\}$ be the subset of $j$ where $l_jn_j\not=0$, equipped with its natural order. Using successively for $j\in J$ either \eqref{ver1} or \eqref{ver2}, according to the sign of $l_j$, we obtain with~:

\begin{equation}
\label{bee}
b:=\sum_{j\in J}b_j\lambda^{\sum_{k\in J,k<j}l_kn_k}\({\frac{1-\lambda^{l_jn_j}}{1-\lambda^{n_j}}}\),
\end{equation}

\smallskip
\noindent
the following relation, for a new constant $C'$ and all $K,L$ large enough~:

$$\sum_{k=K}^{K+L}\|\alpha\lambda^{-k}b+\theta_{k-1}-\theta_k\|^2\leq C'.$$

\smallskip
\noindent
Now, for any $n_j\not=0$, whatever the sign of $n_j$ is, we arrive at, for some constant $C'$ and all $K,L$ large enough~:

$$\sum_{k=K}^{K+L}\|\alpha\lambda^{-k}b\({\frac{1-\lambda^{n_j}}{1-\lambda}}\)+\theta_{k-n_j}-\theta_k\|^2\leq C'.$$

\smallskip
\noindent
Set $b'=b/(1-\lambda)$. Hence, for any $0\leq j\leq N$ with $n_j\not=0$, for some new constant $C'$ and all $K,L$ large enough, using \eqref{pisot}~:

$$\sum_{k=K}^{K+L}\|\alpha\lambda^{-k}(b_j-b'(1-\lambda^{n_j}))\|^2\leq C'.$$

\smallskip
\noindent
Let $0\leq j\leq N$, with $n_j\not=0$. If $b_j\not=b'(1-\lambda^{n_j})$, then we deduce again (still by Cassels \cite{cassels}, chap. 8, Theorems I and II) that $1/\lambda$ is a Pisot number and $b_j=b'(1-\lambda^{n_j})+(1/\alpha)\mu_j$, with $\mu_j\in{\cal T}(1/\lambda)$. The other case is $b_j=b'(1-\lambda^{n_j})$. In any case, we obtain that for all $0\leq j\leq N$~:

\begin{equation}
\label{conj}
\varphi_j(x)=b'+\lambda^{n_j}(x-b')+(1/\alpha)\mu_j,
\end{equation}

\smallskip
\noindent
for some $\mu_j\in{\cal T}(1/\lambda)$. Finally, remark that \eqref{conj} says that the $(\varphi_j)_{0\leq j\leq N}$ are conjugated with the $(\psi_j)_{0\leq j\leq N}$, where $\psi_j(x)=\lambda^{n_j}x+\mu_j$. Precisely $\varphi_j=f\circ \psi_j\circ f^{-1}$, with $f(x)=x/\alpha+b'$.
 
\medskip
\noindent
- {\it Case 2.} Any $\varphi_j$ with $n_j=0$ is the identity. The conclusion is the same, because there now necessarily exists some $0\leq j\leq N$ with $n_j\not=0$ and $b_j\not=b'(1-\lambda^{n_j})$, otherwise $b'$ is a common fixed point for all $(\varphi_j)_{0\leq j\leq N}$.

\medskip
This ends the proof of Theorem \ref{part1}.

\fin

\section{Proof of Theorem \ref{part2}}

Let $N\geq 1$ and affine maps $\varphi_k(x)=\lambda^{n_k}x+\mu_k$, for $0\leq k\leq N$, with $1/\lambda>1$ a Pisot number, relatively prime integers $(n_k)_{0\leq k\leq N}$ and $\mu_k\in {\cal T}(1/\lambda)$, for $0\leq k\leq N$. Let $p\in{\cal C}_N$ and denote by $(\varepsilon_n)_{n\in\Z}$ a two-sided family of $i.i.d.$ random variables with law $p$, to which again the probability $\P$ and the expectation $\E$ refer. We suppose that $\E(n_{\varepsilon_0})>0$. Without loss of generality, $n_N\leq\cdots\leq n_0$ and in particular $n_0\geq1$. For general background on Markov chains, cf Spitzer \cite{spitz}.

\medskip
Recall the cocycle notations for the $(n_{\varepsilon_i})_{i\in\Z}$ introduced before the statement of the theorem and denote by $\theta$ the formal shift such that $\theta \varepsilon_l=\varepsilon_{l+1}$, $l\in\Z$. We have for all $k$ and $l$ in $\Z$~:

$$S_{k+l}=S_k+\theta^kS_l.$$

\smallskip
\noindent
Then $\nu$ is the law of $X=\sum_{l\geq0}\mu_{\varepsilon_l}\lambda^{S_l}$. We write $Q\in\Z[X]$ for the minimal polynomial of $1/\lambda$, of degree $s+1$, with roots $\alpha_0=1/\lambda$, $\alpha_1,\cdots, \alpha_s$, where $|\alpha_k|<1$, for $1\leq k\leq s$. The case $s=0$ corresponds to $1/\lambda$ an integer $\geq2$ (using then usual conventions regarding sums or products). Recall that for any $k\in\Z$, $\sum_{l\in\Z}\mu_{\varepsilon_l}\lambda^{k+S_l}\mod 1$ is a well-defined $\T$-valued random variable.

\bigskip
\noindent
{\it Step 1.} In order to prove the convergence in law of $(\lambda^{-n}X,\lambda^{-n-1}X,\cdots,\lambda^{-n-s}X)\mod \Z^{s+1}$, as $n\rightarrow+\infty$, it is enough to prove, for any $(m_0,\cdots,m_s)\in\Z^{s+1}$, the convergence of~:

$$\E\({e^{2i\pi\sum_{0\leq u\leq s}m_u\lambda^{-n-u}X}}\)=\E\({e^{2i\pi \sum_{l\geq0}(\alpha\mu_{\varepsilon_l})\lambda^{-n+S_l}}}\),$$

\smallskip
\noindent
with $\alpha=\sum_{0\leq u\leq s}m_u\lambda^{-u}$. Notice that $\alpha\mu_j\in{\cal T}(1/\lambda)$, for $0\leq j\leq N$. We make the proof when $\alpha=1$, the one for $\alpha$ being obtained by changing $(\mu_j)$ into $(\alpha\mu_j)$.

\medskip
\noindent
Since $\sum_{l<0}\mu_{\varepsilon_l}\lambda^{-n+S_l}\mod1$ converges a.-s. to 0 in $\T$, as $n\rightarrow+\infty$, it is enough to consider expectations with $\sum_{l\in\Z}\mu_{\varepsilon_l}\lambda^{-n+S_l}\mod1$ in the exponential. Let $k\in\Z$ be a fixed integer. For $n\geq0$, that will tend to $+\infty$, consider $(S_l)_{l\in\Z}$ and the first $q\in\Z$ such that $S_q\geq n$. We have~:

\begin{eqnarray}
\E\({e^{2i\pi \sum_{l\in\Z}\mu_{\varepsilon_l}\lambda^{k-n+S_l}}}\)&=&\sum_{0\leq r<n_0}\sum_{q\in\Z}\E\({e^{2i\pi \sum_{l\in\Z}\mu_{\varepsilon_l}\lambda^{(k-n+S_q)+(S_l-S_q)}}1_{S_{q-u}<n,u\geq1,S_{q}=n+r}}\)\nonumber\\
&=&\sum_{0\leq r<n_0}\sum_{q\in\Z}\E\({e^{2i\pi \sum_{l\in\Z}\mu_{\varepsilon_l}\lambda^{k+r+\theta^qS_{l-q}}}1_{\theta^qS_{-u}<-r,u\geq1,\theta^qS_{-q}=-n-r}}\)\nonumber\\
&=&\sum_{0\leq r<n_0}\sum_{q\in\Z}\E\({e^{2i\pi \sum_{l\in\Z}\mu_{\varepsilon_{l-q}}\lambda^{k+r+S_{l-q}}}1_{S_{-u}<-r,u\geq1,S_{-q}=-n-r}}\)\nonumber\\
&=&\sum_{0\leq r<n_0}\sum_{q\in\Z}\E\({e^{2i\pi \sum_{l\in\Z}\mu_{\varepsilon_{l}}\lambda^{k+r+S_{l}}}1_{S_{-u}<-r,u\geq1,S_{-q}=-n-r}}\).\nonumber
\end{eqnarray}

\medskip
\noindent
For each $0\leq r<n_0$, we now observe that we can move the sum $\sum_{q\in\Z}$ inside the expectation, using the theorem of Fubini, if we first show the finiteness of~:

\begin{eqnarray}
\sum_{q\in\Z}\E\({1_{S_{-q}=-n-r}}\)&=&\E\({\sum_{q\geq0}1_{S_{-q}=-n-r}}\)+\E\({\sum_{q\geq1}1_{S_{q}=-n-r}}\).\nonumber
\end{eqnarray}

\smallskip
\noindent
This is true, since, as soon as $n$ is larger than some constant (because of the missing term for $q=0$ in the second sum), this equals $G^-(0,-n-r)+G^+(0,-n-r)<+\infty$, where $G^-(x,y)$ and $G^+(x,y)$ are the Green functions, finite for every integers $x$ and $y$, respectively associated to the $i.i.d.$ transient random walks $(S_{-q})_{q\geq0}$ and $(S_{q})_{q\geq0}$. Let $\sigma^+_k$, for $k\in\Z$, be the first time $\geq0$ when $(S_{q})_{q\geq0}$ touches $k$. We have $G^+(x,y)=\P_0(\sigma^+_{y-x}<\infty)G^+(0,0)$. With some symmetric quantities, one has $G^-(x,y)=\P_0(\sigma^-_{y-x}<\infty)G^-(0,0)$. 

\medskip
We therefore obtain~:

$$\E\({e^{2i\pi \sum_{l\in\Z}\mu_{\varepsilon_l}\lambda^{k-n+S_l}}}\)=\sum_{0\leq r<n_0}\E\({e^{2i\pi \sum_{l\in\Z}\mu_{\varepsilon_{l}}\lambda^{k+r+S_{l}}}1_{S_{-u}<-r,u\geq1}\({\sum_{q\in\Z}1_{S_{-q}=-n-r}}\)}\).$$

\medskip
\noindent
Let us now fix $0\leq r< n_0$ and consider the corresponding term of the right-hand side. First of all, for $n>0$ larger than some constant (so that $S_0\not=-n-r$)~:

\begin{eqnarray}
\E\({\sum_{q<0}1_{S_{-q}=-n-r}}\)=\P_0(\sigma^+_{-n-r}<\infty)G^+(0,0)\rightarrow0,
\end{eqnarray}

\smallskip
\noindent
as $n\rightarrow+\infty$, since $(S_{q})_{q\geq0}$ is transient to the right. We thus only need to consider~:

$$T(-n):=\E\({e^{2i\pi \sum_{l\in\Z}\mu_{\varepsilon_{l}}\lambda^{k+r+S_{l}}}1_{S_{-u}<-r,u\geq1}N(-n-r)}\),$$

\medskip
\noindent
where $N(-k-r):=\sum_{q\geq0}1_{S_{-q}=-n-r}$. Consider an integer $M_0$, that will tend to $+\infty$ at the end. The difference of $T(-n)$ with the following expression~:

$$\E\({e^{2i\pi \sum_{l\geq-M_0}\mu_{\varepsilon_{l}}\lambda^{k+r+S_{l}}}1_{S_{-u}<-r,1\leq u\leq M_0}N(-n-r)}\)$$

\medskip
\noindent
is bounded by $A+B$, where, first~:

\begin{eqnarray}
A&=&\E\[{\left|{e^{2i\pi \sum_{l\in\Z}\mu_{\varepsilon_l}\lambda^{k+r+S_l}}-e^{2i\pi \sum_{l\geq-M_0}\mu_{\varepsilon_l}\lambda^{k+r+S_l}}}\right|N(-n-r)}\]\nonumber\\
&=&\E\[{\left|{1-e^{2i\pi \sum_{l<-M_0}\mu_{\varepsilon_l}\lambda^{k+r+S_l}}}\right|N(-n-r)}\]\nonumber\\
&\leq &\({\E\[{\left|{1-e^{2i\pi \sum_{l<-M_0}\mu_{\varepsilon_l}\lambda^{k+r+S_l}}}\right|^2}\]}\)^{1/2}\({\E(N(-n-r)^2)}\)^{1/2}\nonumber\\
&\leq &\({\E\[{\left|{1-e^{2i\pi \sum_{l<-M_0}\mu_{\varepsilon_l}\lambda^{k+r+S_l}}}\right|^2}\]}\)^{1/2}\({\E(N(0)^2)}\)^{1/2},\nonumber\end{eqnarray}

\medskip
\noindent
because $N(-n-r)$ is stochastically dominated by $N(0)$. Notice that $N(0)$ is square integrable, as it has exponential tail. The first term on the right-hand side also goes to 0, as $M_0\rightarrow+\infty$, by dominated convergence. The other term $B$ is~:

\begin{eqnarray}
B&=&\E\({1_{S_{-u}<-r,1\leq u\leq M_0,\exists v>M_0, S_{-v}\geq-r}N(-n-r)}\)\nonumber\\
&\leq &\P(\exists v>M_0, S_{-v}\geq-r)^{1/2}\({\E(N(-n-r)^2)}\)^{1/2}\nonumber\\
&\leq &\P(\exists v>M_0, S_{-v}\geq-r)^{1/2}\({\E(N(0)^2)}\)^{1/2},\nonumber\end{eqnarray}

\medskip
\noindent
as before. The first term on the right-hand side goes to 0, as $M_0\rightarrow+\infty$, since $(S_{-v})$ is transient to $-\infty$, as $v\rightarrow+\infty$. As a result~: 

$$T(-n)=\E\({e^{2i\pi \sum_{l\geq-M_0}\mu_{\varepsilon_{l}}\lambda^{k+r+S_{l}}}1_{S_{-u}<-r,1\leq u\leq M_0}N(-n-r)}\)+o_{M_0}(1),$$

\medskip
\noindent
where $o_{M_0}(1)$ goes to 0, as $M_0\rightarrow+\infty$, uniformly in $n$. Now, when $n>0$ is large enough, $N(-k-r)=\sum_{q\geq0}1_{S_{-q}=-n-r}=\sum_{q\geq M_0}1_{S_{-q}=-n-r}$, for all $\omega$. Taking inside the expectation the conditional expectation with respect to the $\sigma$-algebra generated by the $(\varepsilon_l)_{l\geq -M_0}$, we obtain~:

$$T(-n)=\E\({e^{2i\pi \sum_{l\geq-M_0}\mu_{\varepsilon_{l}}\lambda^{k+r+S_{l}}}1_{S_{-u}<-r,1\leq u\leq M_0}G^-(S_{-M_0},-n-r)}\)+o_{M_0}(1).$$

\medskip
\noindent
Now, things are simpler because $G^-(S_{-M_0},-n-r)$ is bounded by the constant $G^-(0,0)$. Hence, for some new $o_{M_0}(1)$, with the same properties~:

$$T(-n)=\E\({e^{2i\pi \sum_{l\in\Z}\mu_{\varepsilon_{l}}\lambda^{k+r+S_{l}}}1_{S_{-u}<-r,u\geq 1}G^-(S_{-M_0},-n-r)}\)+o_{M_0}(1).$$

\medskip
\noindent
Since $G^-(S_{-M_0},-n-r)\rightarrow1/\E(n_{\varepsilon_0})$, as $n\rightarrow\infty$, by renewal theory (since the $(n_j)$ are relatively prime and $p_j>0$, for all $0\leq j\leq N$; cf Woodroofe \cite{wood}, chap. 2, thm 2.1), staying bounded by $G^-(0,0)$, we get by dominated convergence and next $M_0\rightarrow+\infty$~: 

$$\mbox{lim}_{n\rightarrow+\infty}T(-n)=\frac{1}{\E(n_{\varepsilon_0})}\E\({e^{2i\pi \sum_{l\in\Z}\mu_{\varepsilon_{l}}\lambda^{k+r+S_{l}}}1_{S_{-u}<-r,u\geq 1}}\).$$

\medskip
\noindent
From the initial expression, the limit, if existing, had to be independent on the parameter $k$. So this gives the announced convergence and invariance, hence proving item $i)$ in Theorem \ref{part2}.

\bigskip
\noindent
{\it Step 2.} We now consider the proof of Theorem \ref{part2} $ii)$ and suppose that $\nu$ is continuous. We first show that $m$ is a continuous measure. For a continuous $f:\T\rightarrow\R^+$ and any $k\in\R$, we have~:

$$\int_{\T}f(x)~d{m}(x)\leq \frac{1}{\E(n_{\varepsilon_0})}\sum_{0\leq r<n^*}\E\left[{f(Z_{k+r})}\right].$$

\smallskip
\noindent
Letting $k\in\Z$, we have $Z_k=\sum_{l<0}\mu_{\varepsilon_l}\lambda^{k+S_l}+\lambda^{k}X\mod 1$. Since ${\cal L}(\lambda^kX)$ on $\R$ is continuous, ${\cal L}(\lambda^kX\mod 1)$ on $\T$ is continuous. Since $\sum_{l<0}\mu_{\varepsilon_l}\lambda^{k+S_l}\mod1$ and $\lambda^{k}X\mod 1$ are independent random variables, the law of $Z_k$ on $\T$ is continuous. Thus $m$ is a continuous measure (hence $\cal M$). 

\medskip
\noindent
More generally, if $0\not=n=(n_0,\cdots,n_s)^t\in\Z^{s+1}$ and if $Z$ is random variable with law $\cal M$, then the law of $\langle Z,n\rangle$ on $\T$ is $m_{\alpha}$, measure corresponding to $m$ when replacing the $(\mu_j)$ by $(\alpha\mu_j)$, thus the $(\varphi_j)$ by the $(\psi_j)$, with $\psi_j(x)=\lambda^{n_j}x+\alpha\mu_j$, where $\alpha=\sum_{0\leq u\leq s}n_u\lambda^{-u}$. Since $\alpha\not=0$, because $(\lambda^{-u})_{0\leq u\leq s}$ is a basis of $\Q[\lambda]$ over $\Q$, the $(\psi_j)$ do not have a common fixed point and thus $m_{\alpha}$ is continuous, by the previous reasoning.

\medskip
\noindent
Suppose now that the $(\varphi_j)$ have a common fixed point $c$. Hence $\mu_j=c(1-\lambda^{n_j})$, $0\leq j\leq N$, and $\nu=\delta_c$. Necessarily $c\in\Q[\lambda]$, since the $n_j$ are not all zero. We shall show that $\lambda^{-n}c\mod 1$ converges to a rational number in $\T$, as $n\rightarrow+\infty$. First of all, for $n$ large enough, for all $0\leq j\leq N$~:

$$Tr_{1/\lambda}(c\lambda^{-n})-Tr_{1/\lambda}(c\lambda^{-n+n_j})=Tr_{1/\lambda}(\lambda^{-n}\mu_j)\in\Z.$$

\medskip
\noindent
Hence, for any fixed sequence $(k_j)_{0\leq j\leq N}$, for $n$ large enough, for all $0\leq j\leq N$~:

$$Tr_{1/\lambda}(c\lambda^{-n})-Tr_{1/\lambda}(c\lambda^{-n+k_jn_j})\in\Z.$$  

\medskip
\noindent
Supposing that $\sum_{0\leq j\leq N}k_jn_j=1$, using the previous expression successively $n$ replaced by $n,n-k_0n_0,\cdots,n-\sum_{0\leq j\leq N-1}k_jn_j$, respectively with $j=0,j=1,\cdots,j=N$, and finally adding the results, we obtain that for some large $K>0$, for all $n>K$~:

$$Tr_{1/\lambda}(c\lambda^{-n})-Tr_{1/\lambda}(c\lambda^{-n+1})\in\Z.$$ 

\medskip
\noindent
Let $Tr_{1/\lambda}(c\lambda^{-K})=p/q$. For $n>K$, there exists an integer $l_n$ such that $Tr_{1/\lambda}(c\lambda^{-n})=p/q+l_n$. As a result, denoting by $c=c_0,c_1,\cdots,c_s$ the conjugates of $c$ corresponding to $\Q[\lambda]$ (reminding that $(\alpha_j)_{0\leq j\leq s}$ are that of $1/\lambda=\alpha_0$), we get~:

$$c\lambda^{-n}=p/q+l_n-\sum_{1\leq j\leq s}c_j\alpha_j^n.$$

\smallskip
\noindent
Consequently $\lambda^{-n}c\mod 1$ converges to $p/q$ in $\T$, as $n\rightarrow+\infty$, as announced.

\bigskip
\noindent
{\it Step 3.} Consider the proof of Theorem \ref{part2} $iii)$. We show that when $\nu$ is Rajchman, then ${\cal M}={\cal L}_{\T^{s+1}}$. Fix any $0\not=(n_0,\cdots,n_s)^t\in\Z^{s+1}$ and set $\beta=\sum_{0\leq u\leq s}n_u\lambda^{-u}$. Again $\beta\not=0$. We have~:

$$\sum_{0\leq u\leq s}n_u(\lambda^{-n-u}X)=\beta\lambda^{-n}X.$$

\smallskip
\noindent
Since $\nu$ is Rajchman, $\E(e^{2i\pi\beta\lambda^{-n}X})\rightarrow0$, as $n\rightarrow+\infty$. As a result, the Fourier coefficient of ${\cal M}$ corresponding to $(n_0,\cdots,n_s)$ is zero. Hence ${\cal M}={\cal L}_{\T^{s+1}}$. This implies that $m={\cal L}_{\T}$.

\medskip
To complete the proof of $iii)$, we show that $\nu\perp{\cal L}_{\R}$ implies ${\cal M}\perp{\cal L}_{\T^{s+1}}$. Recall that $Z_{k}=\sum_{l\in\Z}\mu_{\varepsilon_l}\lambda^{k+S_l}\mod1$. For any $f\in C(\T^{s+1},\R)$ and $k\in\Z$~:

$$\frac{1}{\E(n_{\varepsilon_0})}\sum_{0\leq r<n^*}\E\left[{f(Z_{-k+r},Z_{-k+r-1},\cdots,Z_{-k+r-s})1_{S_{-v}<-r,v\geq1}}\right]=\int_{\T^{s+1}}f(x)~d{\cal M}(x),$$ 

\smallskip
\noindent
with $n^*=\max_{0\leq j\leq N} n_j$. We now fix $k\geq n^*$ so that $Tr_{1/\lambda}(\lambda^{-l}\mu_j)\in\Z$, for $0\leq j\leq N$, $l\geq k-n^*$. 

\medskip
\noindent
For $0\leq j\leq N$, denote by $(\mu_j^{(t)})_{0\leq t\leq s}$ the conjugates of $\mu_j=\mu_j^{(0)}$ corresponding to the field $\Q[\lambda]$. Let $0\leq r<n^*$. Taking any $0\leq u\leq s$ and $l<0$, we have~:

$$\mu_{\varepsilon_l}\lambda^{-u-k+r+S_l}=Tr_{1/\lambda}(\mu_{\varepsilon_l}\lambda^{-u-k+r+S_l})-\sum_{1\leq t\leq s}\mu^{(t)}_{\varepsilon_l}\alpha_t^{u+k-r-S_l}.$$

\smallskip
\noindent
The role of the indicator function is now fundamental. On the event $\{S_{-v}<-r,v\geq1\}$, we have $Tr_{1/\lambda}(\mu_{\varepsilon_l}\lambda^{-u-k+r+S_l})\in\Z$, by our choice of $k$, since $l\leq-1$. As a result, introducing the {\it real} random variables~:

\begin{equation}
\label{iyr}
Y_u^{(r)}=\lambda^{-u}\sum_{l\geq0}\mu_{\varepsilon_l}\lambda^{-k+r+S_l}-\sum_{1\leq t\leq s}\alpha_t^{u+k-r}\sum_{l<0}\mu^{(t)}_{\varepsilon_l}\alpha_t^{-S_l},
\end{equation}

\smallskip
\noindent
together with $Y^{(r)}=(Y_0^{(r)},\cdots,Y_s^{(r)})$, we obtain that for any $f\in C(\T^{s+1},\R)$~:

\begin{equation}
\label{absm}
\frac{1}{\E(n_{\varepsilon_0})}\sum_{0\leq r<n^*}\E\left[{f(Y^{(r)})1_{S_{-v}<-r,v\geq1}}\right]=\int_{\T^{s+1}}f(x)~d{\cal M}(x).
\end{equation}

\smallskip
\noindent
Hence, for any $f\in C(\T^{s+1},\R^+)$~:

\begin{equation}
\label{abs}
\int_{\T^{s+1}}f(x)~d{\cal M}(x)\leq\frac{1}{\E(n_{\varepsilon_0})}\sum_{0\leq r<n^*}\E\left[{f(Y^{(r)})}\right].
\end{equation} 

\smallskip
\noindent
Fix any $0\leq r<n^*$ and let $X_0=\sum_{l\geq0}\mu_{\varepsilon_l}\lambda^{-k+r+S_l}$ and for $1\leq j\leq s$, $X_j=-\sum_{l<0}\mu^{(j)}_{\varepsilon_l}\alpha_j^{k-r-S_l}$. By definition, $(Y^{(r)})^t=V(X_0,\cdots,X_s)^t$, where $V$ is the Vandermonde matrix~:

$$V=\left({\begin{array}{cccc}
1&1&\cdots&1\\
\lambda^{-1}&\alpha_1&\cdots&\alpha_s\\
\vdots&\vdots&\vdots&\vdots\\
\lambda^{-s}&\alpha_1^s&\cdots&\alpha_s^s\\
\end{array}}\right).$$

\medskip
\noindent
The matrix $V$ is invertible (since the roots of the minimal polynomial $Q$ of $1/\lambda$ are simple). By Cramer's formula~: 

$$X_0=\sum_{0\leq i\leq s}\gamma_iY_i^{(r)},$$

\smallskip
\noindent
with $\gamma_i=\mbox{det}(V^{(i)})/\mbox{det}(V)$, where $V^{(i)}$ is obtained from $V$ by replacing the first column by $e_i$, denoting by $(e_i)_{0\leq i\leq s}$ the canonical basis of $\R^{s+1}$. 

\medskip
\noindent
Notice now that each $\gamma_i$ is real (first of all, $1/\lambda$ is a real root of $Q$; next, regrouping the other roots in conjugate pairs, when conjugating $\gamma_i$ one gets permutations in the numerator $\mbox{det}(V^{(i)})$ and the denominator $\mbox{det}(V)$, the same ones, so $\bar{\gamma_i}=\gamma_i$). As $V$ is invertible, $\gamma:=(\gamma_i)_{0\leq i\leq s}\not=0$.

\medskip
\noindent
We have $X_0=\langle Y^{(r)},\gamma\rangle$. Since $\nu$ is singular with respect to ${\cal L}_{\R}$, we also have ${\cal L}(X_0)\perp {\cal L}_{\R}$, as $X_0=\lambda^{-k+r}X$. As $\gamma\not=0$, we get that ${\cal L}(Y^{(r)})\perp {\cal L}_{\R^{s+1}}$. As a result, ${\cal L}(Y^{(r)}\mod\Z^{s+1})\perp {\cal L}_{\T^{s+1}}$, for all $0\leq r<n^*$. Finally, \eqref{abs} implies that ${\cal M}\perp {\cal L}_{\T^{s+1}}$, as announced.

\medskip
This ends the proof of Theorem \ref{part2}.

\fin

\section{Proof of Theorem \ref{part3}}

The context is the same as that of Theorem \ref{part2}, but now the $(\varphi_k)_{0\leq k\leq N}$ are strict contractions. Precisely, let $N\geq 1$ and $\varphi_k(x)=\lambda^{n_k}x+\mu_k$, for $0\leq k\leq N$, with $1/\lambda>1$ a fixed Pisot number, relatively prime integers $(n_k)_{0\leq k\leq N}$, with now $n_0\geq\cdots \geq n_N\geq1$, without loss of generality, and $\mu_k\in {\cal T}(1/\lambda)$, for $0\leq k\leq N$. 

\bigskip
\noindent
{\it Step 1.} We first show Theorem \ref{part3} $i)$, using again the arguments appearing in the previous section. If $\nu$ is absolutely continuous with respect to ${\cal L}_{\R}$, then ${\cal M}={\cal L}_{\T^{s+1}}$. The event $\{S_{-v}<0,v\geq1\}$ has this time probability one. Looking at \eqref{absm} with $r=0$, we get that the law of $Y^{(0)}\mod\Z^{s+1}$ is absolutely continuous with respect to ${\cal L}_{\T^{s+1}}$, with a density bounded by $\E(n_{\varepsilon_0})$. Hence the law of $Y^{(0)}$ on $\R^{s+1}$ is absolutely continuous with respect to ${\cal L}_{\R^{s+1}}$, with a density also bounded by $\E(n_{\varepsilon_0})$. Since the $(\varphi_k)_{0\leq k\leq N}$ are strict contractions, the $n_j$ are $\geq1$, so the random variable $Y^{(0)}$ is evidently bounded, cf \eqref{iyr}. As a result the density of the law of $Y^{(0)}$ with respect to ${\cal L}_{\R^{s+1}}$ is bounded and with compact support in $\R^{s+1}$. Hence this is also the case of $X_0=\langle Y^{(0)},\gamma\rangle$, where $\gamma=(\gamma_i)_{0\leq i\leq s}\not=0$ is the first line of the inverse of the Vandermonde matrix $V$. Therefore this is also verified for $X=\lambda^{k-r}X_0$. This ends the proof of Theorem \ref{part3} $i)$.

\bigskip
We turn to the proof of Theorem \ref{part3} $ii)$. We shall focus on some Fourier coefficient $\hat{m}(n)$, thus for some fixed $n\in\Z$, of the measure $m$ appearing in Theorem \ref{part2} $i)$. We study its regularity as a function of $p\in{\cal C}_N$, showing its real-analytic character. We then conclude the proof of Theorem \ref{part3} $ii)$ using a theorem on the structure of the set of zeros of a non constant real-analytic function.

\medskip
\noindent
{\it Step 2.} Considering $p\in{\cal C}_N$, denote by $(\varepsilon_n)_{n\in\Z}$ a sequence of $i.i.d.$ random variables with law $p$. Let us fix an integer $n\not=0$, whose exact value will be precised at the end of the proof. We focus on the Fourier coefficient $\hat{m}(n)$ of the measure $m$ introduced in Theorem \ref{part2} $i)$. Let us write $m_p$ in place of $m$ to mark the dependence in $p\in{\cal C}_N$. As $n_j\geq1$, for $0\leq j\leq N$, we have the simplified expression for this Fourier coefficient~:

$$\hat{m}_p(n)=\frac{1}{\E(n_{\varepsilon_0})}\Delta_p\mbox{, with }\Delta_p=\Delta_p(k)=\sum_{0\leq r<n_0}\E\left({e^{2i\pi n\sum_{l\in\Z}\mu_{\varepsilon_l}\lambda^{k+r+S_l}}1_{n_{\varepsilon_{-1}}>r}}\right),$$

\smallskip
\noindent
where this last quantity is independent on $k\in\Z$, by Theorem \ref{part2} $i)$. The expectation $\E(n_{\varepsilon_0})$ also depends on $p$, but to study the zeros of $p\longmapsto \hat{m}_p(n)$ we just need to focus on $\Delta_p$. We now consider the regularity of $p\longmapsto \Delta_p$ on the domain ${\cal C}_N$.

\medskip
For any $k\in\Z$, observe first that $\Delta_p(k)$ is well-defined, with the same formula as above, on the closure $\bar{\cal C}_N$. Fixing $k\in\Z$, the map $p\longmapsto \Delta_p(k)$ is continuous on $\bar{\cal C}_N$, as this function is the uniform limit on $\bar{\cal C}_N$, as $L\rightarrow+\infty$, of the continuous maps~:

$$p\longmapsto \sum_{0\leq r<n_0}\E\left({e^{2i\pi n\sum_{-L\leq l\leq L}\mu_{\varepsilon_l}\lambda^{k+r+S_l}}1_{n_{\varepsilon_{-1}}>r}}\right).$$ 

\smallskip
\noindent
It follows that $p\longmapsto \Delta_p(k)=\Delta_p$ is well-defined on $\bar{\cal C}_N$, continuous and independent on $k$. We shall now prove using standard methods that it is in fact real-analytic in a classical sense, precised below. Let us take $k=0$ and fix $0\leq r<n_0$. Using independence, write~: 

\begin{eqnarray}
\E\left({e^{2i\pi n\sum_{l\in\Z}\mu_{\varepsilon_l}\lambda^{r+S_l}}1_{n_{\varepsilon_{-1}}>r}}\right)&=&\E\left({e^{2i\pi n\sum_{l\geq0}\mu_{\varepsilon_l}\lambda^{r+S_l}}}\right)\E\left({e^{2i\pi n\sum_{l\leq-1}\mu_{\varepsilon_l}\lambda^{r+S_l}}1_{n_{\varepsilon_{-1}}>r}}\right).\nonumber
\end{eqnarray}

\medskip
\noindent
Call $F(p)$ and $G(p)$ respectively the terms appearing in the right-hand side. We shall show that both functions are real-analytic functions of $p$. This property will be inheritated by $p\longmapsto \Delta_p$. We treat the case of $p\longmapsto F(p)$, the case of $G(p)$ needing only to rewrite first the $\mu_{\varepsilon_l}\lambda^{r+S_l}$, appearing in the definition of $G(p)$, as soon as $l<0$ is large enough (depending only the $(\mu_j)_{0\leq j\leq N}$, since $n_k\geq 1$, for all $k$), as $-\sum_{1\leq j\leq s}\alpha_j^{-r-S_l}\mu_{\varepsilon_l}^{(j)}$, quantity equal to $\mu_{\varepsilon_l}\lambda^{r+S_l}$ in $\T$, where the $(\mu_k^{(j)})_{1\leq j\leq s}$ are the conjugates of $\mu_k$  corresponding to the field $\Q[\lambda]$.

\medskip
Fix now $p\in\bar{\cal C}_N$. Let $\N=\{0,1,\cdots\}$ and the symbolic space $S=\{0,\cdots,N\}^{\N}$, equipped with the left shift $\sigma$. For $x=(x_0,x_1,\cdots)\in S$, we define~:

$$g(x)=e^{2i\pi n\({\sum_{l\geq0}\mu_{x_l}\lambda^{r+n_{x_0}+\cdots+n_{x_{l-1}}}}\)}.$$

\smallskip
\noindent
Introducing the product measure $\mu_p=(\sum_{0\leq j\leq N}p_j\delta_j)^{\otimes\N}$ on $S$, we can write~:

$$F(p)=\int_{S}g~d\mu_p.$$

\smallskip
\noindent
Denote by $C(S)$ the space of continuous functions $f:S\rightarrow\C$ and introduce the operator $P_p:C(S)\rightarrow C(S)$ defined by~:

\begin{equation}
\label{pp}
P_p(f)(x)=\sum_{0\leq j\leq N}p_jf((j,x)),~x\in S,
\end{equation}

\smallskip
\noindent
where $(j,x)\in S$ is the word obtained by the left concatenation of the symbol $j$ to $x$. The operator $P_p$ is Markovian, i.e. $f\geq0\Rightarrow P_p(f)\geq0$ and verifies $P_p{\bf 1}={\bf 1}$, where ${\bf 1}(x)=1$, $x\in S$. The measure $\mu_p$ has the invariance property $\int_SP_p(f)~d\mu_p=\int_Sf~d\mu_p$,~$f\in C(S)$. For $f\in C(S)$ and $k\geq0$, introduce the variation~:

$$\mbox{Var}_k(f)=\sup\{|f(x)-f(y)|,~(x,y)\in S^2,x_i=y_i,~0\leq i<k\}.$$

\medskip
\noindent
For any $0<\alpha<1$, let $|f|_{\alpha}=\sup\{\alpha^{-k}\mbox{Var}_k(f),~k\geq0\}$, as well as $\|f\|_{\alpha}=|f|_{\alpha}+\|f\|_{\infty}$. We denote by ${\cal F}_{\alpha}$ the complex Banach space of fonctions $f$ on $S$ such that $\|f\|_{\alpha}<\infty$. Any ${\cal F}_{\alpha}$ is preserved by $P_p$. Observe now that $g\in {\cal F}_{\alpha}$ for $\lambda\leq\alpha<1$. We fix $\alpha=\lambda$. 

\medskip
As a classical fact from Spectral Theory, cf for example Baladi \cite{baladi}, the operator $P_p:{\cal F}_{\lambda}\rightarrow{\cal F}_{\lambda}$ satisfies a Perron-Frobenius theorem. Let us show this elementarily. For $f\in{\cal F}_{\lambda}$, we have~:

$$P^n_pf(x)=\sum_{0\leq j_1,\cdots,j_n\leq N}p_{j_1}\cdots p_{j_n}f((j_1,\cdots,j_n,x)).$$

\smallskip
\noindent
This furnishes $\mbox{Var}_k(P^n_pf-{\bf 1}\int_Sf~d\mu_p)=\mbox{Var}_k(P^n_pf)\leq\mbox{Var}_{k+n}(f)$. Therefore~:

$$\left|{P^n_p(f)-{\bf 1}\int_Sf~d\mu_p}\right|_{\lambda}\leq \lambda^n|f|_{\lambda}.$$

\smallskip
\noindent
In a similar way, we can write~:

\begin{eqnarray}
(P^n_pf-{\bf 1}\int_Sf~d\mu_p)(x)&=&P^n_p(f)(x)-{\bf 1}(x)\int_SP_p^n(f)~d\mu_p\nonumber\\
&=&\sum_{0\leq j_1,\cdots,j_n\leq N}p_{j_1}\cdots p_{j_n}\int_S(f((j_1,\cdots,j_n,x))-f((j_1,\cdots,j_n,y)))~d\mu_p(y).\nonumber
\end{eqnarray}

\noindent
Consequently, $\|P^n_pf-{\bf 1}\int_Sf~d\mu_p\|_{\infty}\leq \mbox{Var}_n(f)\leq\lambda^n|f|_{\lambda}$. Putting things together, finally~:

$$\|P^n_p(f-{\bf 1}\int_Sf~d\mu_p)\|_{\lambda}\leq 2\lambda^n\|f\|_{\lambda}.$$
 
\smallskip
\noindent
This shows that $1$ is a simple eigenvalue and that the rest of the spectrum of $P_p$ is contained in the closed disk of radius $\lambda<1$. Remark that this holds uniformly on $p\in\bar{\cal C}_N$.

\medskip
\noindent
Fix some circle $\Gamma$ centered at 1 and with radius $0<r<1-\lambda$. By standard functional holomorphic calculus, cf Kato \cite{kato}, for any $p\in\bar{\cal C}_N$, the following operator, involving the resolvent, is a continuous (Riesz) projector on $\mbox{Vect}({\bf1})$~:

\begin{equation}
\label{proj}
\Pi_p=\int_{\Gamma}(zI-P_p)^{-1}dz.
\end{equation}

\noindent
Moreover $\Pi_p({\cal F}_{\lambda})$ and $(I-\Pi_p)({\cal F}_{\lambda})$ are closed $P_p$-invariant subspaces, with~:

$${\cal F}_{\lambda}=\Pi_p({\cal F}_{\lambda})\oplus(I-\Pi_p)({\cal F}_{\lambda}).$$

\smallskip
\noindent
Also, in restriction to $(I-\Pi_p)({\cal F}_{\lambda})$, the spectral radius of $P_p$ is less than $\lambda$. 

\medskip
\noindent
Recall that $N\geq1$. We view a function of $p\in \bar{\cal C}_N$ in terms of the first $N$ variables $(p_0,\cdots,p_{N-1})\in\R^N$. Let $\eta'=(\eta_0,\cdots,\eta_{N-1})$ and $\eta=(\eta_0,\cdots,\eta_{N-1},-(\eta_0+\cdots+\eta_{N-1}))$. For any $p\in\bar{\cal C}_N$ and any $\eta'$ (even when $p+\eta\not\in\bar{\cal C}_N$), we can define the continuous operator $P_{p+\eta} : {\cal F}_{\lambda}\rightarrow{\cal F}_{\lambda}$ by \eqref{pp}. It always verifies the relation~:

$$P_{p+\eta}=P_p+\sum_{0\leq j\leq N-1}\eta_jQ_j,$$
  
\noindent
where $Q_j(f)(x)=f(j,x)-f(N,x)$. Denote by $B_N(0,\delta)$ the open Euclidean ball in $\R^N$ of radius $\delta$. Let $\lambda<\lambda'<1-r$. For any $p$ in $\bar{\cal C}_N$, there exists $\delta>0$ such that when $\eta'\in B_N(0,\delta)$, then $\bf1$ is still a simple eigenfunction of $P_{p+\eta}$, with $P_{p+\eta}\bf 1=\bf1$, the rest of the spectrum of $P_{p+\eta}$ being contained in the disk of radius $\lambda'$ and $\Pi_{p+\eta}$, also defined by \eqref{proj}, is a continuous projector on $\mbox{Vect(\bf1})$; this follows from the implicit function theorem, cf Rosenbloom \cite{rosen}, Kato \cite{kato}. By compacity of $\bar{\cal C}_N$, we can choose $\delta>0$ uniformly on $p\in\bar{\cal C}_N$. This defines some open $\delta$-neighborhood ${\cal C}^{\delta}_N$ of $\bar{\cal C}_N$.

\medskip
When $p\in\bar{\cal C}_N$, we have $\int_Sf~d\mu_p=0$, for $f\in(I-\Pi_p)({\cal F}_{\lambda})$. Thus for any $f\in{\cal F}_{\lambda}$~:

$$\Pi_p(f)=\({\int_Sf~d\mu_p}\){\bf1}.$$

\medskip
\noindent
Applying this to the function $g$ of interest to us, we obtain that when $p\in\bar{\cal C}_N$~:

$$F(p){\bf1}=\int_{\Gamma}(zI-P_p)^{-1}(g)dz.$$

\medskip
\noindent
The function $F$ is next extended to ${\cal C}^{\delta}_N$ by the previous formula. Recall the following definition :

\begin{defi} 

$ $

\noindent
A function $h:{\cal C}^{\delta}_N\rightarrow\C$, seen as a function of $(p_0,\cdots,p_{N-1})$, admits a development in series around $p\in{\cal C}^{\delta}_N$, if there exists $\varepsilon>0$ such that for $\eta'=(\eta_0,\cdots,\eta_{N-1})\in B_N(0,\varepsilon)$ and writing $\eta=(\eta',-(\eta_0+\cdots+\eta_{N-1}))$, then $h(p+\eta)$ is given by an absolutely converging series~:

$$h(p+\eta)=\sum_{l_0\geq0,\cdots,l_{N-1}\geq0}A_{l_0,\cdots,l_{N-1}}\eta_0^{l_0}\cdots\eta_{N-1}^{l_{N-1}}.$$

\noindent
A function is real-analytic in ${\cal C}^{\delta}_N$ if it admits a development in series around all $p\in{\cal C}^{\delta}_N$.
\end{defi}

Let us now check that $p\longmapsto F(p)$ is real-analytic on ${\cal C}^{\delta}_N$ in the previous sense. Let $p\in {\cal C}^{\delta}_N$. For $z\in\Gamma$ and $\eta'$ small enough (and the corresponding $\eta$), we can write~:

\begin{eqnarray}
(zI-P_{p+\eta})^{-1}&=&\({I-(zI-P_p)^{-1}\sum_{0\leq j\leq N-1}\eta_jQ_j}\)^{-1}(zI-P_{p})^{-1}\nonumber\\
&=&\sum_{n\geq0}\sum_{0\leq j_1,\cdots,j_n\leq N-1}\eta_{j_1}\cdots\eta_{j_n}(zI-P_{p})^{-1}Q_{j_1}\cdots(zI-P_{p})^{-1}Q_{j_n}(zI-P_{p})^{-1}.\nonumber
\end{eqnarray}

\medskip
\noindent
For small enough $\eta'$, uniformly in $z\in\Gamma$, this is absolutely convergent in the Banach operator algebra. We rewrite it as~:

$$(zI-P_{p+\eta})^{-1}=\sum_{l_0\geq0,\cdots,l_{N-1}\geq0}B_{l_0,\cdots,l_{N-1}}(z)\eta_0^{l_0}\cdots\eta_{N-1}^{l_{N-1}},$$

\noindent
converging for the operator norm, uniformly in $z\in\Gamma$. Hence, for small enough $\eta'$ (and thus $\eta$)~:

\begin{eqnarray}
F(p+\eta){\bf 1}&=&\int_{\Gamma}(zI-P_{p+\eta})^{-1}(g)~dz=\sum_{l_0\geq0,\cdots,l_{N-1}\geq0}\eta_0^{l_0}\cdots\eta_{N-1}^{l_{N-1}}\int_{\Gamma}B_{l_0,\cdots,l_{N-1}}(z)(g)~dz.\nonumber
\end{eqnarray}

\noindent
Applying this equality at some particular $x\in S$, we obtain the desired development in series around $p$. This completes this step.

\bigskip
\noindent
{\it Step 3.} Maybe restricting $\delta>0$, taking into account the finite number of functions appearing in the expression of $\Delta_p$, we obtain that $p\longmapsto\Delta_p$ is real-analytic on ${\cal C}^{\delta}_N$. We shall show that if $n\not=0$ has been appropriately chosen at the beginning, then $\Delta_p$ is not zero at some extremal points of $\bar{\cal C}_N$. The point will be that if ever $\Delta_p$ has a zero on $\bar{\cal C}_N$, then this will imply that either $p\longmapsto\mbox{Re}(\Delta_p)$ or $p\longmapsto\mbox{Im}(\Delta_p)$ is non-constant on ${\cal C}^{\delta}_N$.

\medskip
Now if $h:{\cal C}^{\delta}_N\rightarrow\R$ is real-analytic and non-constant, Lojasiewicz's stratification theorem (cf Krantz-Parks \cite{krantzpark}, theorem 5.2.3) says that the real-analytic set $\{p\in{\cal C}^{\delta}_N~|~h(p)=0\}$ is locally a finite union of real-analytic graphs of dimension $\leq N-1$ (points if $N=1$). By compacity of $\bar{\cal C}_N$, the set $\{p\in\bar{\cal C}_N~|~h(p)=0\}$ is included in a finite union of real-analytic graphs of dimension $\leq N-1$.

\medskip
For the sequel, let us write $x\equiv y$ for equality of $x$ and $y$ in $\T$.

\begin{lemme}

$ $

\noindent
Let $d\geq1$ and $\mu\in{\cal T}(1/\lambda)$. The series $\sum_{l\in\Z}\mu\lambda^{ld}\mod1$, well-defined as an element of $\T$, equals a rational number modulo $1$.
\end{lemme}

\noindent
{\it Proof of the lemma~:}

\noindent
Let $l_0\geq1$ be such that $Tr_{1/\lambda}(\lambda^{-l}\mu)\in\Z$, for $l>l_0$. Denote by $(\mu^{(j)})_{0\leq j\leq s}$ the conjugates of $\mu$, with $\mu^{(0)}=\mu$, and $\alpha_1,\cdots,\alpha_{s}$ that of $\alpha_0=1/\lambda$. We have the following equalities on the torus~:

\begin{eqnarray}
\sum_{l\in\Z}\mu\lambda^{ld}&\equiv&\frac{\mu\lambda^{-l_0d}}{1-\lambda^d}+\sum_{l>l_0}\mu\lambda^{-ld}\equiv\frac{\mu\lambda^{-l_0d}}{1-\lambda^d}-\sum_{1\leq i\leq s}\mu^{(i)}\sum_{l>l_0}\alpha_i^{ld}\equiv\frac{\mu\lambda^{-l_0d}}{1-\lambda^d}-\sum_{1\leq i\leq s}\mu^{(i)}\frac{\alpha_i^{(l_0+1)d}}{1-\alpha_i^d}\nonumber\\
&\equiv&-\({\frac{\mu\lambda^{-(l_0+1)d}}{1-\lambda^{-d}}+\sum_{1\leq i\leq s}\mu^{(i)}\frac{\alpha_i^{(l_0+1)d}}{1-\alpha_i^d}}\)=-Tr_{1/\lambda}\({\frac{\mu\lambda^{-(l_0+1)d}}{1-\lambda^{-d}}}\)\in\Q.
\nonumber\end{eqnarray}

\fin

\bigskip
We complete the argument. Fixing $0\leq j\leq N$ and $p^{j}=(0,\cdots,0,1,0,\cdots,0)$, where the 1 is at place $j$, we have for $k\in\Z$, recalling that $1\leq n_j\leq n_0$~:

\begin{eqnarray}
\Delta_{p^j}=\Delta_{p^j}(k)&=&\sum_{0\leq r< n_0}e^{2i\pi n\sum_{l\in\Z}\mu_j\lambda^{k+r+ln_j}}1_{n_j>r}=\sum_{0\leq r< n_j}e^{2i\pi n\sum_{l\in\Z}\mu_j\lambda^{k+r+ln_j}}.\nonumber
\end{eqnarray}

\smallskip
\noindent
Notice in passing that the invariance with respect to $k$ is now obvious, as we sum over $r$ on a full period of length $n_j$. Now, taking $k=0$, we have~:

$$\Delta_{p^j}=\sum_{0\leq r< n_j}e^{2i\pi n(A_{j,r}/B_{j,r})},$$

\smallskip
\noindent
for rational numbers $A_{j,r}/B_{j,r}$, making use of the previous lemma, since $\lambda^r\mu_j\in{\cal T}(1/\lambda)$, for any $r$. If for example $n$ is a multiple of $B_{j,r}$ for any $0\leq r<n_j$, we get $\Delta_{p^j}=n_j\geq1$, which gives what was desired. This ends the proof of the theorem.

\fin

\bigskip
\noindent
\begin{remark}
Lojasiewicz's stratification theorem, giving the local structure of $\{p\in{\cal C}^{\delta}_N~|~h(p)=0\}$, is a difficult theorem. In an elementary way, using the implicit function theorem, one can show that the set of zeros of a real-valued real analytic non constant function is locally included in a countable union of connected real-analytic graphs of codimension one. 
\end{remark}

\bigskip
\noindent
\begin{remark}
In the general case, when the $(\varphi_k)_{0\leq k\leq N}$ are not all strict contractions, the method seems to reach some limit. Using the notation ${\cal D}_N(r)$ of the Introduction, with $r=(\lambda^{n_k})_{0\leq k\leq N}$, and considering as in {\it Step 2} the regularity of $p\longmapsto F(p)$ on ${\cal D}_N(r)$, it is not difficult to show continuity, using some standard coupling argument. The real-analytic character, if ever true, a priori requires more work. Still setting $S=\{0,\cdots,N\}^{\N}$ and $\mu_p=(\sum_{0\leq j\leq N}p_j\delta_j)^{\otimes\N}$ on $S$, we again have~:

$$F(p)=\int_{S}g~d\mu_p,$$

\noindent
with $g(x)=e^{2i\pi n\({\sum_{l\geq0}\mu_{x_l}\lambda^{r+n_{x_0}+\cdots+n_{x_{l-1}}}}\)}$, but this function is only defined $\mu_p$-almost-everywhere. 

\end{remark}

\section{Complements}

\subsection{A numerical example}
Considering an example as simple as possible which is not homogeneous, take $N=1$ and the two contractions $\varphi_0(x)=\lambda x,~\varphi_1(x)=\lambda^2x+1$, where $1/\lambda>1$ is a Pisot number, with probability vector $p=(p_0,p_1)$. Then $n_0=1$, $n_1=2$ and $\nu$ is the law of $\sum_{l\geq0}\varepsilon_l\lambda^{n_{\varepsilon_0}+\cdots+n_{\varepsilon_{l-1}}}$, with $(\varepsilon_n)_{n\geq0}$ $i.i.d.$, with common law $\mbox{Ber}(p_1)$, i.e. $\P(\varepsilon_0=1)=p_1$ and $\P(\varepsilon_0=0)=1-p_1$. We shall take $0\leq p_1\leq 1$ as parameter for simulations. Notice that $\E(n_{\varepsilon_0})=p_0+2p_1=1+p_1$, 

\medskip
Taking $n=1$, $k\in\Z$ and $r\in\{0,1\}$, let us define~:

$$F_{p}(k)=\E\({e^{2i\pi\lambda^k\sum_{l\geq0}\varepsilon_l\lambda^{n_{\varepsilon_0}+\cdots+n_{\varepsilon_{l-1}}}}}\),~G_p(k,r)=\E\({e^{2i\pi\sum_{l\geq0}\varepsilon_l\lambda^{k-(n_{\varepsilon_0}+\cdots+n_{\varepsilon_l})}}1_{n_{\varepsilon_0}>r}}\),$$

\smallskip
\noindent
leading to $\Delta_p=F_p(k)G_p(k,0)+F_p(k+1)G_p(k+1,1)$, for all $k\in\Z$. Writing $m_p$ in place of $m$ for the measure on $\T$ in Theorem \ref{part2} $i)$ (defined when $0<p_1<1$), we get $\hat{m}_p(1)=\Delta_p/(1+p_1)$.  Let us first discuss the choice of probability vector $p=(1-p_1,p_1)$ and Pisot number $1/\lambda$.

\medskip
A degenerated example (the invariant measure being automatically singular with respect to ${\cal L}_{\R}$) is for instance given by $\lambda=(3-\sqrt{5})/2<1/2$. Nevertheless, it is interesting to notice that $\lambda^{-n}\equiv-\lambda^n$, $n\geq0$. Taking $p_1=1/2$, one can check that $\Delta_p=|F_p(1)|^2+|F_p(2)|^2/2$. Necessarily $\Delta_p>0$. Indeed, $k\longmapsto F_p(k)$ verifying a linear recurrence of order two, the equality $\Delta_p=0$ would give $F_p(k)=0$ for all $k$, but $F_p(k)\rightarrow1$, as $k\rightarrow+\infty$. Notice that $(3-\sqrt{5})/2$ is the largest $\lambda$ with this property (it has to be a root of some $X^2-aX+1$, for some integer $a\geq0$). Mention that in general $\Delta_p$ is not real; cf the pictures below.

\medskip
To study an interesting example, we take into account the similarity dimension $s(p,r)$, rewritten here as $s(p,\lambda)$~:

$$s(p,\lambda):=\frac{(1-p_1)\ln (1-p_1)+p_1\ln p_1}{(1-p_1)\ln \lambda+p_1\ln(\lambda^2)}.$$

\medskip
\noindent
The condition $s(p,\lambda)\geq1$ is equivalent to $(1-p_1)\ln (1-p_1)+p_1\ln p_1-(1+p_1)\ln\lambda\leq0$. As a function of $p_1$, the left-hand side has a minimum value $-\ln(\lambda+\lambda^2)$, attained at $p_1=\lambda/(1+\lambda)$. As a first attempt, taking for $1/\lambda$ the golden mean $(\sqrt{5}+1)/2=1,618...$ is in fact not interesting, as in this case $\lambda+\lambda^2=1$, giving $s(p,\lambda)\leq1$. 

\medskip
We instead take (as considered in Section 2) for $1/\lambda$ the Plastic number, i.e. the unique real root of $X^3-X-1$. Approximately, $1/\lambda=1.324718...$. For this $\lambda$~:

$$s(p,\lambda)>1\Longleftrightarrow 0,203...<p_0<0,907....$$

\medskip
\noindent
The other roots of $X^3-X-1=0$ are conjugate numbers $\rho e^{\pm i\theta}$. From the relations $1/\lambda+2\rho\cos\theta=0$ and $(1/\lambda)\rho^2=1$, we deduce $\rho=\sqrt{\lambda}$ and $\cos\theta=-1/(2\lambda^{3/2})$, thus $\theta=\pm2.43...$ rad. For computations, the relations $\lambda^{-n}+\rho^ne^{in\theta}+\rho^ne^{-in\theta}\in\Z$, $n\geq0$, furnish $\lambda^{-n}\equiv-2(\sqrt{\lambda})^n\cos(n\theta)$.

\medskip
\noindent
Let us finally compute the extreme values of $p_1\longmapsto \hat{m}_p(1)$, abusively written as $\hat{m}_{(1,0)}(1)$ and $\hat{m}_{(0,1)}(1)$, since $m_p$ has only been defined for $0<p_1<1$. We first observe that $\hat{m}_{(1,0)}(1)=\Delta_{(1,0)}=F_{(1,0)}(0)G_{(1,0)}(0,0)=1$. At the other extremity~:

\begin{eqnarray}
\Delta_{(0,1)}&=&F_{(0,1)}(0)G_{(0,1)}(0,0)+F_{(0,1)}(1)G_{(0,1)}(1,1)\nonumber\\
&=&e^{2i\pi \sum_{l\geq0}\lambda^{2l}}e^{2i\pi\sum_{l\geq0}\lambda^{-2(l+1)}}+e^{2i\pi\lambda\sum_{l\geq0}\lambda^{2l}}e^{2i\pi\sum_{l\geq0}\lambda^{1-2(l+1)}}\nonumber\\
&=&e^{2i\pi\({\frac{1}{1-\lambda^2}-2\sum_{l\geq0}(\sqrt\lambda)^{2l}\cos(2l\theta)}\)}+e^{2i\pi\({\frac{\lambda}{1-\lambda^2}-2\sum_{l\geq0}(\sqrt\lambda)^{2l+1}\cos((2l+1)\theta)}\)}\nonumber\\
&=&e^{2i\pi\({\frac{1}{1-\lambda^2}-2\mbox{Re}\({\frac{\lambda e^{2i\theta}}{1-\lambda e^{2i\theta}}}\)}\)}+e^{2i\pi\({\frac{\lambda}{1-\lambda^2}-2\mbox{Re}\({\frac{\sqrt{\lambda}e^{i\theta}}{1-\lambda e^{2i\theta}}}\)}\)}.\nonumber
\end{eqnarray}

\medskip
\noindent
A not difficult computation, shortened by the observation that $(1-\lambda e^{2i\theta})(1-\lambda e^{-2i\theta})=1/\lambda$, shows that the arguments in the exponential terms (after the $2i\pi$) are respectively equal to $3$ and $0$, leading to $\Delta_{(0,1)}=2$ and therefore $\hat{m}_{(0,1)}(1)=1$.

\medskip
Recalling that $p=(1-p_1,p_1)$, below are respectively drawn the real-analytic maps $p_1\longmapsto\mbox{Re}(\hat{m}_{p}(1))$, $p_1\longmapsto\mbox{Im}(\hat{m}_{p}(1))$ and the parametric curve $p_1\longmapsto \hat{m}_p(1)$, $0\leq p_1\leq1$.

\medskip
{\centerline{\includegraphics[scale=0.30]{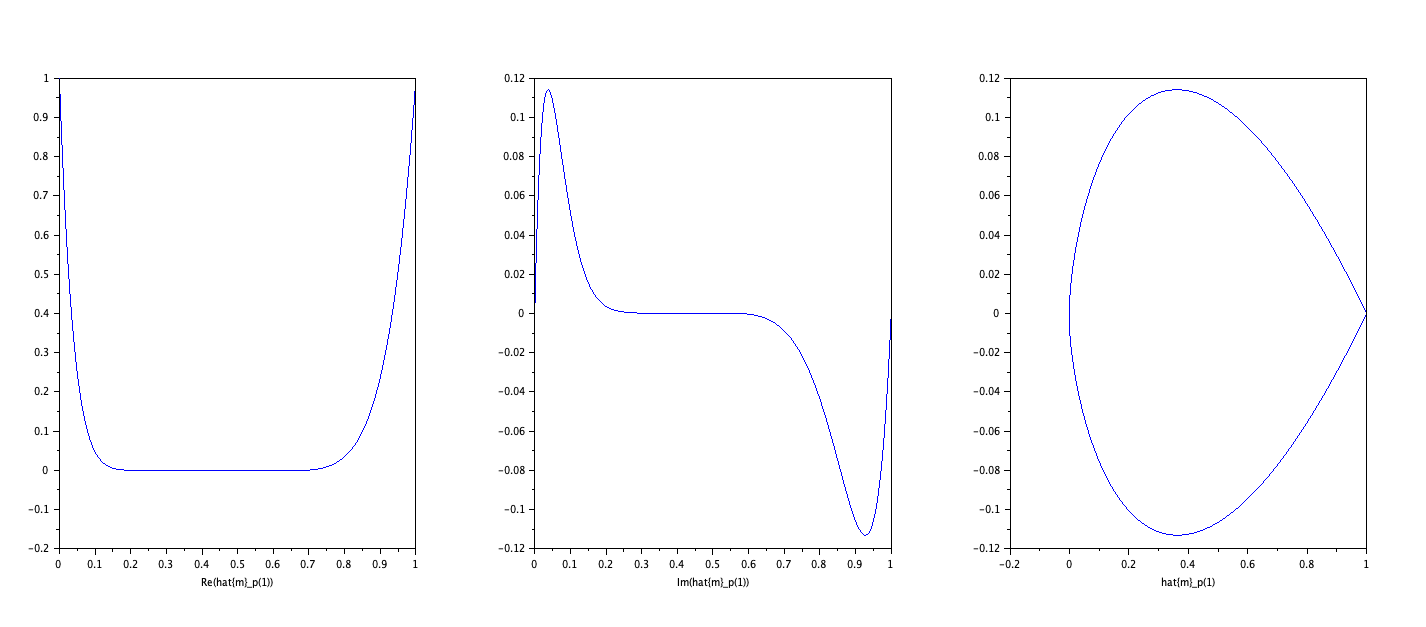}}}

\medskip
\noindent
The first two pictures indicate that $p_1\longmapsto \hat{m}_p(1)$ spends a rather long time near $0$, with $\mbox{Re}(\hat{m}_p(1))$ and $\mbox{Im}(\hat{m}_p(1))$ both around $10^{-4}$.  Let us precise here that one can exploit the product form (given by the exponential) inside the expectation appearing in $F_p(k)$ and $G_p(k,r)$. Using a binomial tree, we make a deterministic numerical computation of $\hat{m}_{p}(1)$, with nearly an arbitrary precision. For example, one can obtain the rather remarquable value :

$$\hat{m}_{(1/2,1/2)}(1)=0,0001186...+i0,0000327...,$$

\medskip
\noindent
where all digits are exact. In this case, $s((1/2,1/2),\lambda)=1,64...>1$. The above pictures were drawn with 1000 points, each one determined with a sufficient precision. This allows to safely zoom on the neighbourhood of $0$ of $p_1\longmapsto \hat{m}_p(1)$, the interesting region. We obtain the following surprising pictures, the one on the right-hand side containing around 500 points :

{\center{\includegraphics[scale=0.30]{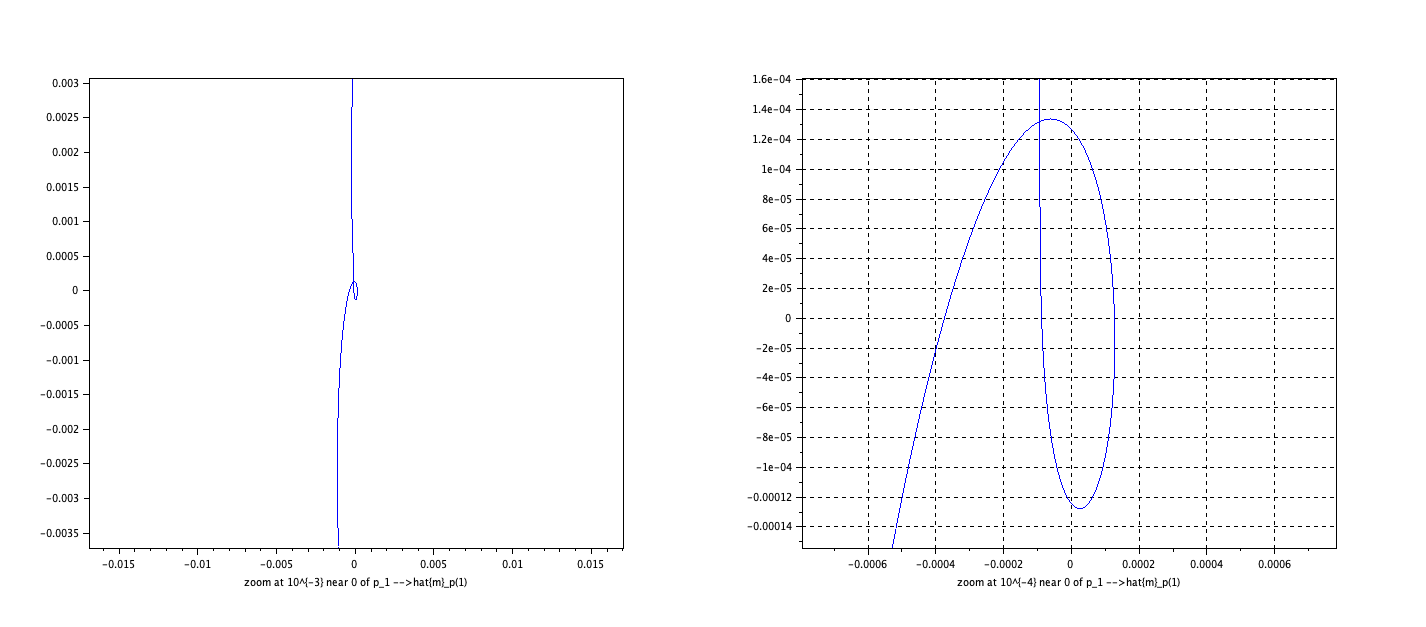}}}

\medskip
\noindent
There are probably profound reasons behind these pictures, that would in particular clarify the condition of non-nullity of the Fourier coefficient $\hat{m}_p(1)$ and more generally of $\hat{m}_p(n)$, $n\in\Z$. Further investigations are necessary, but we can conclude that the curve $p_1\longmapsto\hat{m}_{p}(1)$ is rather convincingly not touching $0$. It may certainly be possible to build a rigorous numerical proof of this fact, but this is not the purpose of the present paper. We informally state~:

\begin{nthm}
$ $

\noindent
Let $N=1$, $0<\lambda<1$, with $1/\lambda>1$ the Plastic number, and $\varphi_0(x)=\lambda x$, $\varphi_1(x)=\lambda^2x+1$. Then for all $p\in{\cal C}_1$, the invariant measure $\nu$ is continuous singular and not Rajchman.

\end{nthm}

\noindent
\begin{remark}
For the same system, but taking for $1/\lambda$ the supergolden ratio, i.e. the fourth Pisot number (the real root of $X^3-X^2-1$), one essentially gets the same pictures.
\end{remark}

\medskip
Still taking for $1/\lambda$ the Plastic number, but for the system $\varphi_0(x)=\lambda^2x$ and $\varphi_1(x)=\lambda^3x+1$, already mentioned in Section 2, recall that the invariant measure $\nu$ is continuous singular and not Rajchman for all $p\in{\cal C}_1$, except when $p=(\lambda^2,\lambda^3)$, in which case $\nu=\frac{1}{1+\lambda}{\cal L}_{[0,1+\lambda]}$. We have drawn below the real analytic curve $p_1\longmapsto\hat{m}_{p}(1)$, with next a zoom at $10^{-3}$ near the origin. This is also interesting, since this time the curve is not self-intersecting, being almost linear near zero and passing at zero exactly for the sole parameter $p_1=\lambda^3$.

{\center{\includegraphics[scale=0.37]{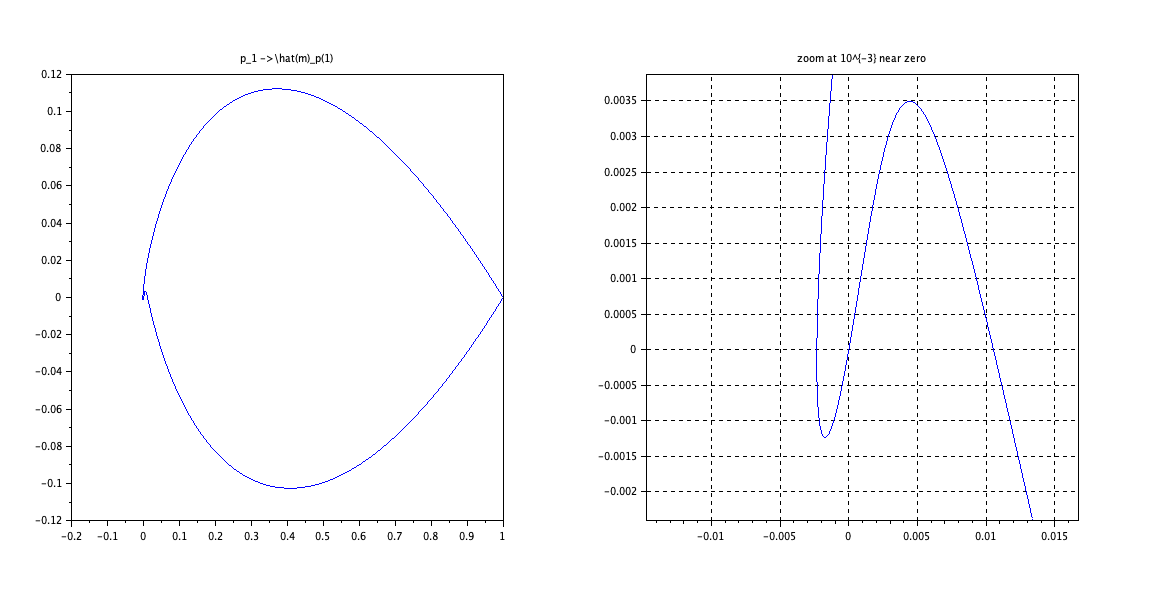}}}

\subsection{Applications to sets of uniqueness for trigonometric series}

Let $N\geq1$ and for $0\leq k\leq N$ affine contractions $\varphi_k(x)=r_kx+b_k$, with reals $(r_k)$ and $(b_k)$ such that $0<r_k<1$ for all $k$. As a general fact, Theorem \ref{part1} has some consequences in terms of sets of multiplicity for trigonometric series, cf for example Salem \cite{salem2} or Zygmund \cite{zyg} for details. As in the Introduction, let $F\subset\R$ be the unique non-empty compact set, verifying the self-similarity relation $F=\cup_{0\leq k\leq N}\varphi_k(F)$. With $\N=\{0,1,\cdots\}$ and $S=\{0,\cdots,N\}^{\N}$, one has~:

\begin{eqnarray}
F&=&\left\{{\sum_{l\geq0}b_{x_l}r_{x_0}\cdots r_{x_{l-1}},~(x_0,x_1,\cdots)\in S}\right\}.\nonumber
\end{eqnarray}

\smallskip
Let us place on the torus $\T$ and consider trigonometric series. Recall that a subset $E$ of $\T$ is a set of uniqueness ($U$-set), if whenever a trigonometric series $\sum_{n\geq 0}(a_n\cos(2\pi x)+b_n\sin(2\pi x))$, with complex numbers $(a_n)$ and $(b_n)$, converges to $0$ for all $x\not\in E$, then $a_n=b_n=0$ for all $n\geq0$. Otherwise $E$ is said of multiplicity ($M$-set). 

\begin{thm}
\label{mult}

$ $

\noindent
Let $N\geq1$ and for $0\leq k\leq N$ affine contractions $\varphi_k(x)=r_kx+b_k$, where $0<r_k<1$, with no common fixed point. Suppose that the system $(\varphi_k)_{0\leq k\leq N}$ is not affinely conjugated to a family in Pisot form. Then $F\mod 1\subset\T$ is a $M$-set. 

\end{thm}

\medskip
\noindent
{\it Proof of the theorem~:}

\noindent
Any $p\in{\cal C}_N$ gives a Rajchman invariant probability measure $\nu$ supported by $F\subset\R$. Hence $F\mod(1)\subset\T$ supports the probability $\tilde{\nu}$, image of $\nu$ under the projection $x\longmapsto x\mod1$, from $\R$ to $\T$. Then $\tilde{\nu}$ is a Rajchman measure on $\T$, so, cf Salem \cite{salem2} (chap. V), $F\mod1$ is a $M$-set. \fin

\medskip
In the other direction, in general more delicate, we shall simply apply existing results. For the following statement, fixing $0<\lambda<1$ and integers $n_k\geq 1$, for $0\leq k\leq N$, notice that for any $(x_0,x_1,\cdots)\in S$, we have $\sum_{l\geq0}\lambda^{n_{x_0}+\cdots+n_{x_{l-1}}}(1-\lambda^{n_{x_l}})=1$.

\begin{thm}
\label{uniq}

$ $

\noindent
Let $N\geq1$ and suppose that the $(\varphi_k)$ are affine contractions of the form $\varphi_k(x)=\lambda^{n_k}x+b_k$, with $b_k=ba_k+c(1-\lambda^{n_k})$, for some $0<\lambda<1$ with $1/\lambda$ a Pisot number $>N+2$, relatively prime positive integers $n_k\geq1$, $0\leq a_k\in\Q[\lambda]$ and real numbers $b\geq0$ and $c$. Then the non-empty compact self-similar set $F=\cup_{0\leq k\leq N}\varphi_k(F)\subset\R$ can be written as $F=bG+c$, where $G$ is the compact set :

$$G=\left\{{\sum_{l\geq0}a_{x_l}\lambda^{n_{x_0}+\cdots+n_{x_{l-1}}},~(x_0,x_1,\cdots)\in S}\right\}.$$

\medskip
\noindent
Assume that $bG\subset[0,1)$, so that $bG$ and $F$ can be seen as subsets of $\T$. Then $F$ is $U$-set.

\end{thm}

\medskip
\noindent
{\it Proof of the theorem~:}

\noindent
Up to replacing $b$ and the $(a_k)$ respectively by $br$ and $(a_k/r)$, for some $r>1$ in $\Q$, we may assume that $0\leq a_k<1/(1-\lambda)$, for all $0\leq k\leq N$. Then~:

$$G\subset H:=\left\{{\sum_{l\geq0}\eta_l\lambda^l,~\eta_l\in\{0,a_0,\cdots,a_N\},~l\geq0}\right\}\subset[0,1).$$

\noindent
Since $1/\lambda>N+2$ is a Pisot number and all $a_0,\cdots,a_N$ are in $\Q[\lambda]$, it follows from the Salem-Zygmund theorem, cf Salem \cite{salem2}, chap. VII, paragraph 3, on perfect homogeneous sets, that $H$ is a perfect $U$-set. Mention that in this theorem, one also assumes that $\max_{0\leq k\leq N}a_k=1/(1-\lambda)$ and that successive $a_u<a_v$ in $[0,1)$ verify $a_v-a_u\geq \lambda$. These conditions serve to give a geometrical description of the perfect homogeneous set $H$ in terms of dissection, without overlaps. They are in fact not used in the proof, where only the above description of $H$ is important (one can indeed start reading Salem \cite{salem2}, chap. VII, paragraph 3, directly from line 9 of the proof).

\medskip
As a subset of a $U$-set, $G$ is also a $U$-set. This is also the case of $bG$, by hypothesis a subset of $[0,1)$, using Zygmund, Vol. I, chap. IX, Theorem 6.18 (the proof, not obvious, is in Vol. II, chap. XVI, 10.25, and relies on Fourier integrals). Hence, $F=bG+c$ is also a $U$-set, as any translate on $\T$ of a $U$-set is a $U$-set. This ends the proof of the theorem.

\fin

\bigskip
\noindent
\begin{remark} As a general fact, the hypothesis $1/\lambda>N+2$ ensures that $H$ and $F$ have zero Lebesgue measure, which is a necessary condition for a set to be a $U$-set. \end{remark}

\bigskip
\noindent
{\it Acknowledgments.} We thank B. Kloeckner for initial discussions on this topic, J. Printems for his very precious help in the numerical part of the last section and K. Conrad for several enlightening discussions concerning Algebraic Number Theory questions. We also are grateful to P. Varju for pointing out an error in a first version of Theorem \ref{part3}. His work with H. Yu \cite{varjuhu}, completes section 6.2 and fully characterizes sets of uniqueness among self-similar sets.

\providecommand{\bysame}{\leavevmode\hbox to3em{\hrulefill}\thinspace}

\bigskip
{\small{\sc{Univ Paris Est Creteil, CNRS, LAMA, F-94010 Creteil, France\\
Univ Gustave Eiffel, LAMA, F-77447 Marne-la-Vall\'ee, France
 }}}

\it{E-mail address~:} {\sf julien.bremont@u-pec.fr}

\end{document}